\documentclass{article}
\usepackage{graphicx}
\usepackage{amssymb}
\usepackage{amsmath}
\usepackage{mathtools}
\usepackage{a4wide}
\usepackage{caption}
\usepackage{subcaption}
\usepackage[affil-it]{authblk}
\usepackage{xr}
\usepackage{color}
\usepackage[utf8]{inputenc}
\usepackage[T1]{fontenc}

\begin{document}

\title{A note on time-optimal paths on perturbed spheroid}
\author{Piotr Kopacz}
\affil{\small{Jagiellonian University, Faculty of Mathematics and Computer Science\\ 6, Prof. S. Łojasiewicza, 30 - 348, Kraków,
Poland}}
\affil{Gdynia Maritime University, Faculty of Navigation \\3,  Al. Jana Pawła II, 81-345, Gdynia, Poland}
\date{\small{\textit{E-mail address:} \texttt{piotr.kopacz@im.uj.edu.pl}}}
\maketitle
\begin{abstract}
\noindent
We consider the Zermelo navigation problem on the ellipsoid of revolution (spheroid) in the presence of a perturbation $W$ determined by a mild velocity vector field, $|W|<1$, with application of Finsler metric of Randers type in the context of the corresponding optimal control represented by a time-efficient ship's heading $\varphi(t)$ (steering direction). As the example we present the solutions to the problem on an oblate ellipsoid with acting infinitesimal rotation. 
\end{abstract}

\

\smallskip
\noindent
\textbf{M.S.C. 2010}: 53C22, 53C60, 53C21, 53B20, 49N90, 49J15, 34H05. 
\smallskip


\smallskip

%%%%%%%%%%%%%%%%%%%%%%%%%%%%%%%%%%%%%%%%%%

\section{Introduction}

A spheroid is widely applied as a geometric model in navigation as well as in geodesy, cartography and engineering. Due to the combined effects of gravity and rotation the Earth's shape is often approximated by an oblate spheroid. In particular, cartographic and geodetic systems for the Earth are based on a reference ellipsoid, for instance the current World Geodetic System. Satellite data have provided new measurements to define the best Earth-fitting ellipsoid and for relating existing coordinate systems to the Earth’s centre of mass. It is the approximate shape of many planets and celestial bodies and the quickly-spinning stars. Several moons of the Solar System approximate prolate spheroids in shape, though they are actually triaxial ellipsoids. The most common shapes for the density distribution of protons and neutrons in an atomic nucleus are spherical and spheroidal, where the polar axis is assumed to be the spin axis or direction of the spin angular momentum vector. Also, the current algorithms of the accurate and global navigational calculations correspond to the spheroidal geometric models (cf. \cite{earle, pallikaris, kuo_guo}). 

The time-optimal paths in the presence of perturbing vector field $W$ on a spheroid is of our interest.  Thus, we refer to the Zermelo navigation problem; see \cite{colleen_shen, chern_shen} for more details. Our goal is to consider the problem on the ellipsoid by means of Finsler geometry in the context of the optimal control when the key role plays the formula for navigating ship's heading (steering direction), i.e. the angle $\varphi=\varphi(t)$ which the vector of the relative velocity forms with a fixed direction. In dimension two and three, with the Euclidean background metric the implicit solutions for the optimal headings were given by Ernst Zermelo (1931) when the problem was formalized initially \cite{zermelo2, zermelo}. Then, in particular, the problem was followed in detailed analysis with application of the Hamiltonian formalism by Constantin Caratheodory \cite{caratheodory}.  
%($h_{ij}):=(\delta_{ij})$ where $\delta_{ij}$ denotes the Kronecker delta
The analysis of the problem simplifies if we observe that it suffices to find the locally optimal solution on the tangent space. We aim to find the deviation of the Riemannian geodesics on an ellipsoid and the corresponding steering angles (optimal control) such that the ship following the resulting paths completes her journey in the least time, under the action of distributed wind $W$.  %In fact, we investigate the solutions to the generalized special case of the originally stated  navigation problem. 

%%%%%%%%%%%%%%%%%%%%%%%%%%%%%%%%%%%%%%%%%%

\section{Unperturbed Riemannian geodesics of a spheroid}

In general, let a pair $ (M,h) $ be a Riemannian manifold where $h = h_{ij}dx^i\otimes dx^j$ is a Riemannian metric and the corresponding norm-squared of tangent vectors $\mathbf{y} \in T_x M$ is denoted by $\left|\mathbf{y} \right|_h^2 = h_{ij}y^iy^j = h(\mathbf{y}, \mathbf{y})$. The unit tangent sphere in each $T_x M$ consists of all tangent vectors $\mathbf{u}$ such that $h(\mathbf{u}, \mathbf{u}) = 1$. We assume that a ship proceeds with the constant speed relative to a perturbation, $|\mathbf{u}|_h=1$, and the problem of time-efficient navigation is treated for the case of a mild wind $W$, namely $|W|_h < 1$ everywhere on $M$. The time-optimal paths on Riemannian manifold in unperturbed scenario are represented by the geodesics of the corresponding Riemannian metric. 
Let $\Sigma^2$ be an ellipsoid embedded in the Euclidean space $\mathbb{R}^3$ referred to Cartesian coordinates $(x^i)$, with axes $2r$, $2r$, $2ar$. The parametrization of $\Sigma^2$ in the spherical coordinate system $(\rho, \phi, \theta)$ which we apply in the paper yields $x=r\sin\theta\cos \phi, y=r\sin\theta \sin \phi, z= ar\cos \theta$, where  the azimuth $\phi\in [0, 2\pi)$ and the inclination $\theta \in [0, \pi]$. The parameter $a>0$ determines the shape of an ellipsoid and as a consequence  the flow of the geodesics on a spheroid which can be oblate ($0<a<1$) or prolate ($a>1$). We equip an ellipsoid with the Riemannian metric induced by the Euclidean metric of $\mathbb{R}^3$. Computing the square distance $ds^2$ between two points $(x^i)$ and $(x^i+dx^i)$ of $\Sigma^2$ expressed in the terms of the spherical coordinates $(\phi, \theta)$ on $\Sigma^2$ one obtains $ds^2=r^2\sin^2\theta (d\phi)^2+r^2(\cos^2\theta+a^2 \sin^2\theta)(d\theta)^2$. Hence, the background Riemannian metric $h=h_{ij}dx^idx^j$ is given by $h_{11}=r^2 \sin ^2\theta, h_{22}=r^2 (\cos^2\theta + a^2\sin^2\theta)$, $h_{12}=h_{21}=0$. Computing nonzero Christoffel symbols of the first and second kind yields $\Gamma _{121} = r^2 \sin \theta  \cos \theta $,  $\Gamma _{211}= r^2\sin \theta \cos \theta $, $\Gamma _{222} =\left(a^2-1\right) r^2 \sin \theta \cos \theta $ and $\Gamma ^1_{12}=\cot\theta$, $\Gamma ^2_{11}=-\frac{\sin \theta  \cos \theta }{a^2 \sin ^2\theta +\cos ^2\theta }$, $\Gamma ^2_{22} = \frac{\left(a^2-1\right) \sin \theta  \cos \theta }{a^2 \sin ^2\theta +\cos ^2\theta }$. We thus get the geodesic equations of $\Sigma^2$ in the following form 

\begin{equation}
\ddot{\phi}+2 \dot{\theta} \dot{\phi} \cot\theta=0,  
\label{sph_geo1}
\end{equation}

%\begin{equation}
%\ddot{\theta}+\frac{\left[\left(a^2-1\right) \dot{\theta}^2 - \dot{\phi}^2\right]\sin \theta \cos \theta}{a^2 \sin ^2\theta+\cos ^2\theta} =0
%\end{equation}

\begin{equation}
(a^2 \sin ^2\theta+\cos ^2\theta)\ddot{\theta}+\sin \theta \cos \theta\left[\left(a^2-1\right) \dot{\theta}^2 - \dot{\phi}^2\right]=0,
\label{sph_geo2}
\end{equation}

\noindent
where the dots indicate derivatives with respect to $t$. From now on we assume that $r:=1$ so the ellipsoid has the semiaxes $(1, 1, a)$. Without loss of generality we shall consider the oblate spheroid with $a:=\frac{3}{4}$ in the presented example. 

\begin{figure}[h]
        \centering
%~\includegraphics[width=0.32\textwidth]{img/sf_geo5}
~\includegraphics[width=0.3\textwidth]{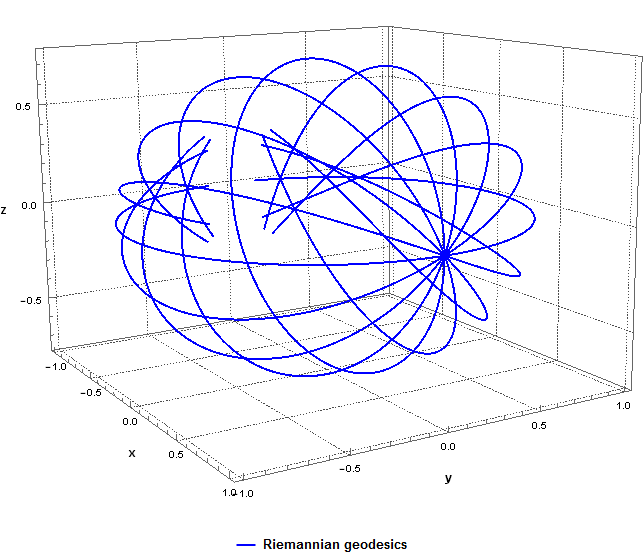}
~\includegraphics[width=0.3\textwidth]{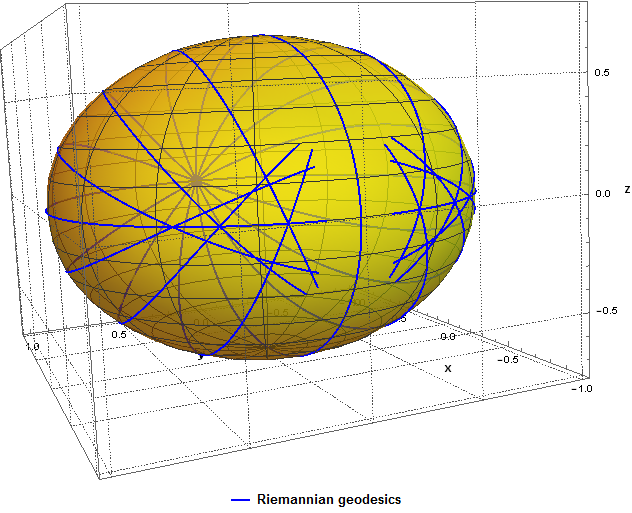}
~\includegraphics[width=0.36\textwidth]{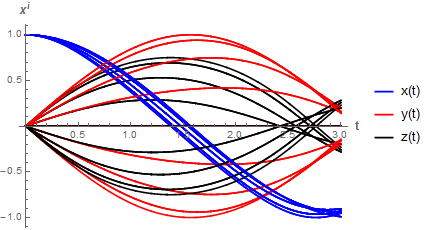}
        \caption{The Riemannian geodesics on the spheroid with $a:=\frac{3}{4}$ starting at $(0, \frac{\pi}{2})$ with the increments $\Delta \varphi_0 = \frac{\pi}{8}$ and their solutions $x(t), y(t), z(t)$ for $(1, 0, 0)\in\mathbb{R}^3$; $t=3$.}
\label{Riem_geo} 
\end{figure}
\noindent
Keeping in mind the assumption that own speed of the ship is unit, the form of the initial conditions in the case of the time-optimal solutions before perturbation become $\phi(0)=\phi_0\in[0, 2\pi)$, $\theta(0)=\theta_0\in[0, \pi]$ and 
\begin{equation}
\dot{\phi}(0)=u, \quad \dot{\theta}(0)=\pm\sqrt\frac{1-u^2\sin^2\theta_0}{\cos^2\theta_0+a^2\sin^2\theta_0}
\end{equation}
or
\begin{equation}
\dot{\phi}(0)=\pm\frac{1}{\sin\theta_0}\sqrt{1-v^2(\cos^2\theta_0+a^2\sin^2\theta_0)}, \quad \dot{\theta}(0)=v. 
\end{equation}
For instance, if a ship commences the passage from any point of the spheroid's equator, i.e. $\theta_0=\frac{\pi}{2}$ then the condition which determines the coordinates of an arbitrary tangent vector is simplified to $u^2+a^2v^2=1$. Let the starting point on $\Sigma^2$ be determined by $\phi_0=0$ and $\theta_0=\frac{\pi}{2}$ to show the flow of unperturbed Riemannian geodesics. Note that, due to the form of set spherical coordinates, if $\dot{\phi}(t)<0$ then a path runs west (clockwise in a top view), and if $\dot{\theta}(t)<0$ then the path runs north. Similarly, for the positive values it goes east and south, respectively. In Figure \ref{Riem_geo} we show the Riemannian geodesics on the spheroid with $a:=\frac{3}{4}$ starting at $(0, \frac{\pi}{2})$, with the increments $\Delta \varphi_0 = \frac{\pi}{8}$ and their solutions $x(t), y(t), z(t)$ for $(1, 0, 0)\in\mathbb{R}^3$, $t=3$.

%\begin{figure}
   %     \centering
%~\includegraphics[width=0.35\textwidth]{img/sf_geo4}
%~\includegraphics[width=0.35\textwidth]{img/sf_geo4a}

   %     \caption{Family - Riemannian geodesic, t=7, $\Delta \varphi_0 = \frac{\pi}{4}$} 
%\end{figure}

%%%%%%%%%%%%%%%%%%%%%%%%%%%%%%%%%%%%%%%%

\section{Introducing perturbing vector field}

Generally, the perturbation applied in the navigation problem can be given in the following form 
\begin{equation}
\label{general_field}
W = \frac{\partial}{\partial t} + W^i\left(t,x^j\right)\frac{\partial}{\partial x^i}
\end{equation}
which includes the time dependence. The wind may have a rotational effect as well as a translational effect on the ship and this effect depends on the ship's heading. Our aim is to consider the stationary wind and this is a special case of  Zermelo's problem. We introduce perturbing vector field in the general form $\tilde{W}(x, y, z)=\tilde{W}^1(x, y, z)\frac{\partial}{\partial x}+\tilde{W}^2(x, y, z)\frac{\partial}{\partial y}+\tilde{W}^3(x, y, z)\frac{\partial}{\partial z}$ which acts on embedded spheroid in $\mathbb{R}^3$. In the new base $(\phi, \theta)$ the form of $W$ becomes $W(\phi, \theta)=W^1(\phi, \theta)\frac{\partial}{\partial\phi}+W^2(\phi, \theta)\frac{\partial}{\partial\theta}$. The condition $|W|_h<1$, where 
\begin{equation}
|W(\phi, \theta)|_h=\sqrt{(W^1(\phi, \theta))^2\sin^2\theta+(W^2(\phi, \theta))^2(\cos^2\theta+a^2\sin^2\theta)}, 
\end{equation}
ensures that $F$ is strongly convex. This is the necessary codition as we study the problem via the construction of the Finsler metric \cite{colleen_shen,  chern_shen}. Let us apply the perturbation $\tilde{W}=cy\frac{\partial}{\partial x}-cx\frac{\partial}{\partial y}+0\frac{\partial}{\partial z}$ which determines a circulation of an angular speed $c$, where $c$ is a constant. Thus, $W$ refers to rigid motions of the spheroid and describes rotation about the $z$-axis. Let $|c|<1$. Then $W$ is given by 
\begin{equation}
	W(\phi, \theta)=-c\frac{\partial}{\partial\phi} \qquad \text{with} \qquad  |W(\phi, \theta)|_h=|c|\sin\theta <1.
\label{ir}
\end{equation}
Consequently, $|W(\phi, \theta)|_h\in[0, c]$ and $|W(\phi,\frac{\pi}{2})|_{max}=c$ in the case of perturbing east wind ($c>0$) running along the ellipsoid's equator.

\begin{figure}[h]
       \centering
~\includegraphics[width=0.3\textwidth]{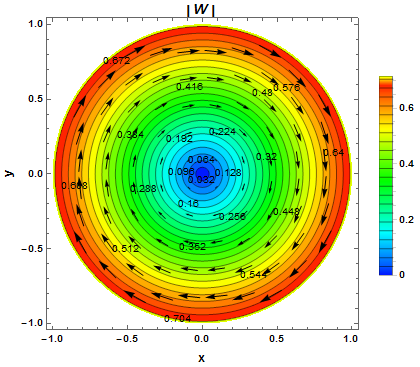}\qquad \quad
~\includegraphics[width=0.4\textwidth]{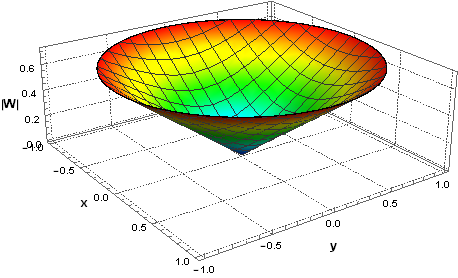}
\caption{The contour plot and the graph of the norm $|W|_h$ in the case of the infinitesimal rotation with $c:=\frac{5}{7}$ as the acting perturbation.} 
\label{ir_fig}
\end{figure}
\noindent
Without loss of generality let us consider east wind (clockwise when looking from the north pole of a spheroid). In the example which follows we shall take $c:=\frac{5}{7}$. Figure $\ref{ir_fig}$ shows the contour plot and the graph of the norm $|W|_h$ of the acting perturbation, namely the infinitesimal rotation. What about the case where $|W|_h\geq1$ ? Note that originally the solutions to the problem with the Euclidean background were not restricted by the wind's strength and also included the cases where  $|W|>|u|$. In the generalizations of the Zermelo navigation problem in Finsler geometry for stronger perturbation, i.e. $|W|_h=1$ the Kropina metric has been applied recently (cf. \cite{kropina}). The strong wind with $|W|_h>1$ has been considered very recently (cf. \cite{sanchez}) in reference to the Lorentzian metric with the general relativity background.

%%%%%%%%%%%%%%%%%%%%%%%%%%%%%%%%%%%%%%%%%%%%%

\section{Time-optimal paths under acting perturbation}

Under the influence of a pertubation $W$ which is modeled by acting vector field, the paths of the shortest time are no longer the geodesics of the Riemannian metric $h$. Instead, they are the geodesics of a special Finsler type $F$, namely Randers metric \cite{colleen_shen, chern_shen}. We shall apply the metric $F$ which correlates to the background geodesics on an ellipsoid. 

\subsection{The resulting metric}

Randers metrics may be identified with the solutions to the navigation problem on Riemannian manifolds. Then a bijection between Randers spaces and pairs $(h,W)$ of Riemannian metrics $h=h_{ij}y^iy^j$ and acting vector fields $W$ on the manifold $M$ is established. 
The resulting Randers metric is composed of the new Riemannian metric and $1$-form, and is given by \cite{colleen_shen}
\begin{equation}
\label{ran}
F(\mathbf{y}) = \frac{ \sqrt{ \left[ h(W,\mathbf{y}) \right]^2 + |\mathbf{y}|^2 \lambda} } {\lambda} - \frac{ h(W,\mathbf{y})} {\lambda}
\end{equation}
where $W_i=h_{ij}W^j$ and $\lambda=1-W^iW_i$. 
Hence, 
\begin{equation}
F(x,\mathbf{y}) = \frac{ \sqrt{( h_{ij}W^iy^j)^2 + |\mathbf{y}|^2 (1-|W|^2)} -h_{ij}W^iy^j} {1-|W|^2}.
\end{equation}

%METRYKA F W 3D Z EUKLID WYJSCIOWA
%\begin{equation} 
%\begin{split}
%F(x,y,z; \mu,\nu, \vartheta) = \frac{-\sqrt{(u W^1+v W^2+w W^3)^2-\left(u^2+v^2+w^2\right) \left((W^1)^2+(W^2)^2+(W^3)^2-1\right)}}{1-|W|_{euklid}^2} \\
%+ \frac{u W^1+v W^2+w W^3}{1-|W|_{euklid}^2}
%\end{split}
%\end{equation}
%where $W^i=W^i(x,y,z)$. 

%\begin{equation}
%\tiny{
%\begin{split}
%F(x,y,z; u,v,w) = \frac{-\sqrt{u^2 \left(-\left((W^2)^2+(W^3)^2-1\right)\right)+2 u W^1 (v W^2+w W^3)-v^2 \left((W^1)^2+(W^3)^2-1\right)+2 v w W^2 W^3-w^2 \left((W^1)^2+(W^2)^2-1\right)}}{1-|W|_{euklid}^2}\\
%+\frac{u W^1+v W^2+w W^3}{1-|W|_{euklid}^2}
%\end{split}
%}
%\end{equation}
\noindent
Being in dimension two after the coordinates' transformation we denote the position coordinates $(x^1, x^2)$ by $(\phi, \theta)$, and expand arbitrary tangent vectors $y^1\frac{\partial}{\partial x^1}+y^2\frac{\partial}{\partial x^2}$ at  $(x^1, x^2)$ as $(\phi, \theta;u,v)$ or $u\frac{\partial}{\partial \phi} + v\frac{\partial}{\partial \theta}$. Thus, adopting the notations in two dimensional case yields

\begin{equation}
 \label{W1W2}
F(\phi, \theta; u,v) = \frac{\sqrt{u^2h_{11}+v^2h_{22}-(uW^2-vW^1)^2h_{11}h_{22}}-uW^1h_{11}-vW^2h_{22}}{1-|W|^2}
 \end{equation}
where $W^i=W^i(\phi, \theta)$. %Może powinno być $y=(u, v, w)$ w R3 a $y=(mi, ni)$ dla sferoidy we wspol sfefycznych, gdbysmy takze zahaczali o R3
The norm of function $F$ mesures travel time on $\Sigma^2$. In our problem the perturbed background Riemannian metric $h$ is the induced Euclidean metric on a spheroid $\Sigma^2$. From \eqref{W1W2} we obtain the Randers metric $F=\alpha+\beta$ on the spheroid $(1, 1, a)$ where  
\begin{equation}
\label{sph1}
%\footnotesize{
\begin{split}
\alpha(\phi, \theta; u, v)=\frac{\sqrt{u^2\sin^2\theta+(\cos^2\theta+(a\sin\theta)^2)[v^2-\sin^2\theta(uW^2(\phi, \theta)-vW^1(\phi, \theta))^2]}}{1-(W^1(\phi, \theta)\sin\theta)^2-(W^2(\phi, \theta))^2[\cos^2\theta+(a\sin\theta)^2]},  
\end{split}
%}
\end{equation}

\begin{equation}
\label{sph2}
%\footnotesize{
\begin{split}
\beta(\phi, \theta; u, v)=-\frac{uW^1(\phi, \theta)\sin^2\theta+vW^2(\phi, \theta)(\cos^2\theta+(a\sin\theta)^2)}{1-(W^1(\phi, \theta)\sin\theta)^2-(W^2(\phi, \theta))^2[\cos^2\theta+(a\sin\theta)^2]}. 
\end{split}
%}
\end{equation}

\

\noindent
For the perturbation $W(\phi, \theta)=-c\frac{\partial}{\partial\phi}$ the form of the resulting metric becomes 
\begin{equation}
\label{sph2d_final}
%\footnotesize{
F(\phi ,\theta; u,v)=\frac{\sqrt{-c^2 v^2 \sin ^2\theta  \left(a^2 \sin ^2\theta +\cos ^2\theta \right)+v^2 \left(a^2 \sin ^2\theta +\cos ^2\theta \right)+u^2 \sin ^2\theta }+c u \sin ^2\theta }{1-c^2 \sin ^2\theta }.
%}
\end{equation}

%\begin{equation}
%\footnotesize{
%L(\phi ,\theta ,u,v)=\frac{\left(\sqrt{-c^2 v^2 \sin ^2(\theta ) \left(a^2 \sin ^2(\theta )+\cos ^2(\theta )\right)+v^2 \left(a^2 \sin ^2(\theta )+\cos ^2(\theta )\right)+u^2 \sin ^2(\theta )}+c u \sin ^2(\theta )\right)^2}{2 \left(1-c^2 \sin ^2(\theta )\right)^2}
%}
%\end{equation}

\noindent
Next, we use the partial derivatives of $L\text{=}\frac{1}{2}F^2$ to obtain the spray coefficients what comes below and consequently the final geodesic equations. Their solutions will determine the time-optimal paths.  

%%%%%%%%%%%%%%%%%%%%%%%%%%%%%%%%%%%%%%%%%%%

\subsection{Time-optimal paths' equations}

A spray on $M$ is a smooth vector field on $TM_0=TM\backslash \{0\}$ locally expressed in the standard  form
\begin{equation}
    \label{spray}
G=y^i\frac{\partial}{\partial x^i}-2G^i \frac{\partial}{\partial y^i},
\end{equation}
where $G^i=G^i(x,\mathbf{y})$ are the local functions on $TM_0$ satisfying $G^i(x, \widehat{\lambda} \mathbf{y})=\widehat{\lambda} ^2G^i(x,\mathbf{y})$, where $\widehat{\lambda} >0$. 
The spray is induced by $F$ and the spray coefficients $G^i$ of $G$ given by \cite{chern_shen}
\begin{equation}
    \label{spray2}
G^i=\frac{1}{4}g^{il}\{ [F^2]_{x^ky^l}y^k-[F^2]_{}x^l\}%=\frac{1}{4}g^{il}\left(2\frac{\partial g_{jl}}{\partial x^k}-\frac{\partial g_{jk}}{\partial x^l}\right)y^jy^k,
\end{equation}
are the spray coefficients of $F$.
In two dimensional case we express the spray coefficients $G^1:=G(\phi,\theta;u,v)$ and $G^2:=H(\phi,\theta;u,v)$ after adopting the spherical coordinates by
\begin{equation}
G(\phi,\theta;u,v)=\frac{\left(\frac{\partial ^2L}{\partial v^2} \frac{\partial L}{\partial \phi}-\frac{\partial L}{\partial \theta} \frac{\partial ^2L}{\partial u\, \partial v}\right)-\frac{\partial L}{\partial v} \left(\frac{\partial ^2L}{\partial \phi\, \partial v}-\frac{\partial ^2L}{\partial \theta\, \partial u}\right)}{2 \left[\frac{\partial ^2L}{\partial u^2} \frac{\partial ^2L}{\partial v^2}-\left(\frac{\partial ^2L}{\partial u\, \partial v}\right)^2\right]},
\end{equation}
\begin{equation}
H(\phi,\theta;u,v)=\frac{\left(\frac{\partial ^2L}{\partial u^2} \frac{\partial L}{\partial \theta}-\frac{\partial L}{\partial \phi} \frac{\partial ^2L}{\partial u\, \partial v}\right)+\frac{\partial L}{\partial u} \left(\frac{\partial ^2L}{\partial \phi\, \partial v}-\frac{\partial ^2L}{\partial \theta\, \partial u}\right)}{2 \left[\frac{\partial ^2L}{\partial u^2} \frac{\partial ^2L}{\partial v^2}-\left(\frac{\partial ^2L}{\partial u\, \partial v}\right)^2\right]}
\end{equation}
Next, we apply the last two formulae to compute the spray coefficients of the resulting Randers metric. However first, we obtain the spray coefficients $G_{\alpha}, H_{\alpha}$ of the new Riemannian term $\alpha$ of $F$ for acting infinitesimal rotation given by \eqref{ir}. The result is 
\begin{equation}
G_{\alpha}(\phi ,\theta; u,v)=\frac{u v  \left(c^2+\csc ^2\theta \right)\cot \theta}{\csc ^2\theta -c^2},  
\end{equation}
% A DLA OGÓLNEGO WIATRU (W1, W2) SPRAYE G, H SĄ POLICZONE, ALE NIECO PRZYDŁUGIE POSTACIE 
%\begin{equation}
%\small{
%H_{\alpha}(\phi ,\theta; u,v)=\frac{\cot \theta  \csc ^4\theta  \left(c^2 \cos 2 \theta  \left(v^2 \left(c^2+a^2-1\right)+u^2\right)-\left(c^2-2\right) v^2 \left(c^2+a^2-1\right)-\left(c^2+2\right) u^2\right)}{4 \left(c^2-\csc ^2\theta \right)^2 \left(a^2+\cot ^2\theta \right)}.
%}
%\end{equation}

\begin{equation}
H_{\alpha}(\phi ,\theta; u,v)=
\frac{ \splitfrac{ \cot \theta  \csc ^4\theta  \left(c^2 \cos 2 \theta  \left(v^2 \left(c^2+a^2-1\right)+u^2\right) \right.}{\left.-\left(c^2-2\right) v^2 \left(c^2+a^2-1\right)-\left(c^2+2\right) u^2\right)}}{4 \left(c^2-\csc ^2\theta \right)^2 \left(a^2+\cot ^2\theta \right)}.
\end{equation}

\begin{figure}[h]
        \centering
~\includegraphics[width=0.3\textwidth]{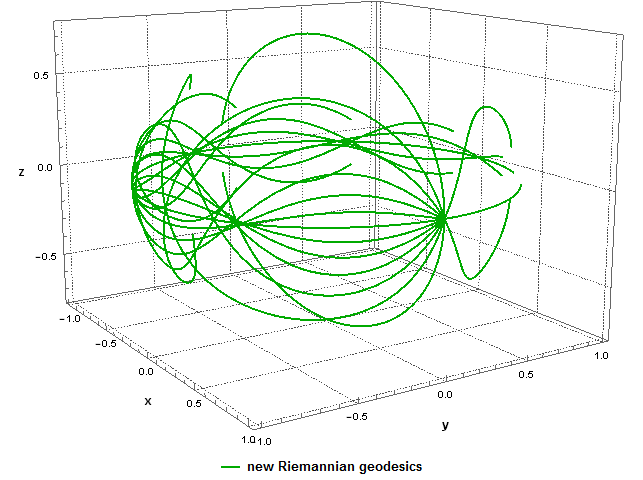}
%~\includegraphics[width=0.35\textwidth]{img/sf_fam_t=3a}
%~\includegraphics[width=0.32\textwidth]{img/sf_alfa1}
~\includegraphics[width=0.32\textwidth]{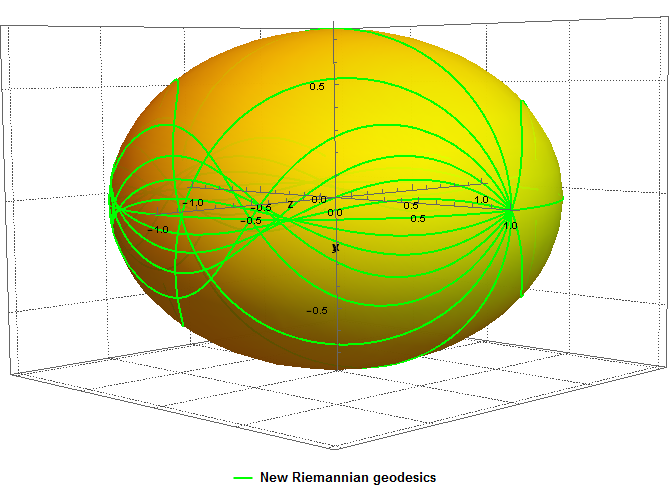}
%~\includegraphics[width=0.3\textwidth]{img/sf_alfa3}

        \caption{The geodesics of the new Riemannian metric $\alpha$ starting from $(0, \frac{\pi}{2})$ with the increments $\Delta \varphi_0=\frac{\pi}{8}$, $t=3$.} 
\label{geo_riem_new}
\end{figure}
\noindent
For a standard local coordinate system $(x^i, y^i)$ in $TM_0$ the geodesic equation for Finsler metric is expressed in the general form
\begin{equation}
\label{geo}
\dot{y}^i+2G^i(x,\mathbf{y})=0. %\qquad
%\ddot{x}^i+2G^i(x,\dot{x})=0
\end{equation}
%Hence, 
%\begin{equation}
%\label{geo1}
%\ddot{x}^i+\frac{1}{2} g^{il} (2\frac{\partial g_{jl}}{\partial x^k} - \frac{\partial g_{jk}}{\partial x^l})\dot{x}^j\dot{x}^k=0.
%\end{equation}

\noindent
Hence, the geodesic equations for the new Riemannian metric $\alpha$ become 

\begin{equation}
\ddot{\phi}+\frac{2 \dot{\theta}  \dot{\phi} \left(c^2+\csc ^2\theta\right)\cot \theta}{\csc ^2(\theta)-c^2}=0,  
\label{alpha_1ex}
\end{equation}

\begin{equation}
\ddot{\theta}+\frac{\cot \theta  \csc ^4\theta  \left(\left(c^2+a^2-1\right) \left(c^2 \cos 2 \theta -c^2+2\right) \dot{\theta}^2+\left(c^2 \cos 2 \theta -c^2-2\right) \dot{\phi}^2\right)}{2 \left(c^2-\csc ^2\theta \right)^2 \left(a^2+\cot ^2\theta \right)}=0.
\label{alpha_2ex}
\end{equation}

\noindent
The geodesics of the new Riemannian metric $\alpha$ starting from $(0, \frac{\pi}{2})$ with the increments $\Delta \varphi_0=\frac{\pi}{8}$ and $t=3$ are presented in Figure \ref{geo_riem_new}.
The solution curves in the problem are found by working out the geodesics of $F$ given by \eqref{sph1} and \eqref{sph2}. Analyses involving Randers spaces are generally difficult and finding solutions to the geodesic equations is not straightforward \cite{brody, chern_shen}. As it takes a while if one computes some quotients manually, we create some programmes with use of Wolfram Mathematica ver. 10.3 to generate the graphs and provide some numeric computations when the complete symbolic ones cannot be obtained. The numerical schemes can give useful information studying the geometric properties of obtained solutions as is shown in the attached graphs. Due to complexity of the obtained spray coefficients $G$ and $H$ of the resulting metric \eqref{sph2d_final} and consequently the final time-optimal paths' equations, we now present the spray coefficients which induce the Randers geodesic equations for $W(\phi, \theta)=-c\frac{\partial}{\partial\phi}$. Let us abbreviate
\begin{equation*}
\psi = \sqrt{\sin ^2\theta  \left(-c^2 a^2 v^2 \sin ^2\theta +a^2 v^2+u^2\right)+v^2 \cos ^2\theta  \left(1-c^2 \sin ^2\theta \right)},
\end{equation*}
\begin{equation*}
\mu = a^2 v^2+u^2,
\end{equation*}
\begin{equation*}
\tau = \left(3 c u^2   \psi  -c a^2 v^2  \psi +3u\mu\right).
\end{equation*}
%koniec symboli
Hence,
\begin{equation}
G(\phi ,\theta; u,v)=
	\frac{
		 v \cos \theta  \left(c 
	\psi 
	+u\right) \left(\csc ^3\theta 
	 \psi 
	+c u \csc \theta \right)^3}
	{
	\splitfrac{ \left(\csc ^2\theta -c^2\right)
	\left(c^3 u \left(u^2-3 a^2 v^2\right)+\csc ^4\theta  \mu 
	\psi 
	 -\frac{1}{8} v^2 \cot ^2\theta  \csc ^4\theta  \right.} {\left.   \left(c^2\cos 2 \theta  
	   -c^2+2\right) 	 
	 \left(-4
	  \psi 
	  +6 c u \cos 2 \theta -6 c u\right)+c \csc ^2\theta  
	  \tau
	  \right)  
	 }
	 },
\end{equation}

\begin{equation}
H(\phi ,\theta; u,v)=-
\frac{
\splitfrac{
2 \sin \theta  \cos \theta  \left(
\psi
+c u \sin ^2\theta \right)^3 \left(c^4 \left((2 a^2-1\right) v^2 \sin ^4 \theta -c^2 \sin ^2\theta  
 \right.}{\left.
  \left(\left(3 a^2-2\right) v^2+u^2\right)-2 c u 
\psi
+c^2 v^2 \left(c^2 \sin ^2\theta -1\right)\cos ^2\theta +a^2 v^2-u^2-v^2\right)
}
}
{
\splitfrac{
\left(c^2 \sin ^2 \theta -1\right)^2  \left(a^2 \sin ^2\theta +\cos ^2\theta \right) \left(v^2 \cos ^2\theta  \left(c^2 \cos 2 \theta -c^2+2\right) \left(-2 
\psi
\right.\right.}{\left.\left.
+3 c u \cos 2 \theta -3 c u\right)
-4 \sin ^2\theta  \left(c^3 u  \left(u^2-3 a^2 v^2\right) \sin ^4\theta +
\mu \psi 
+c \tau \sin ^2 \theta 
	\right)
	\right)
}
}.
\end{equation}
\begin{figure}[h]
        \centering
~\includegraphics[width=0.3\textwidth]{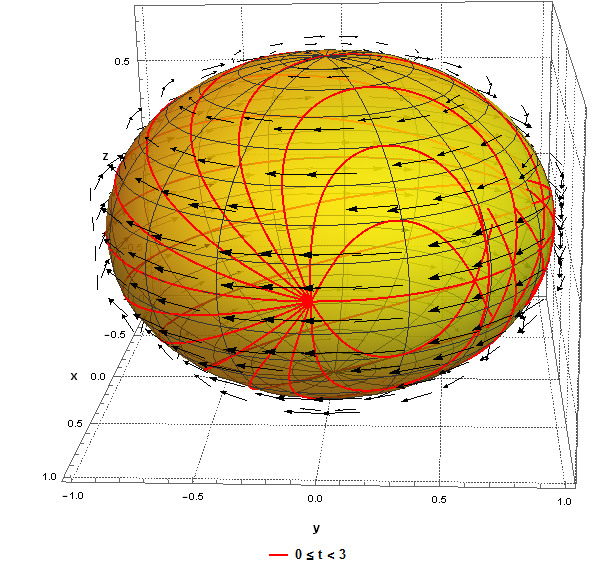}
~\includegraphics[width=0.32\textwidth]{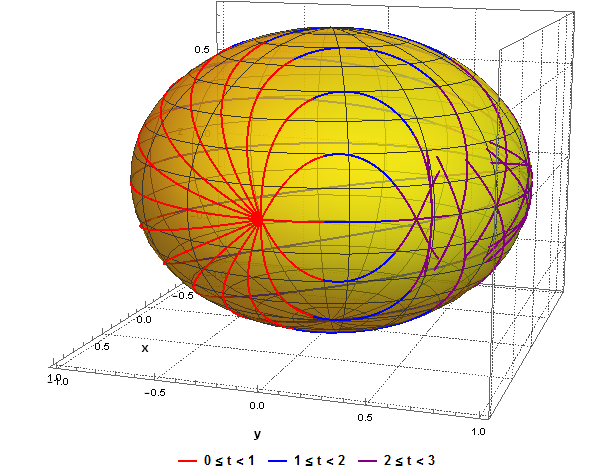}
~\includegraphics[width=0.32\textwidth]{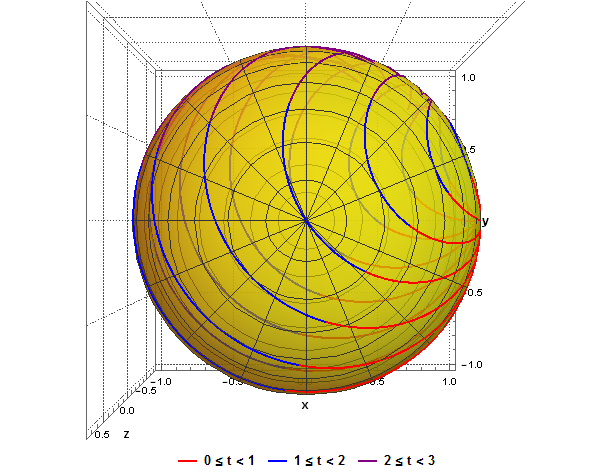}
%~\includegraphics[width=0.35\textwidth]{img/sf_fam_t=3b}
%~\includegraphics[width=0.35\textwidth]{img/sf_fam_t=3}

        \caption{The time-efficient paths on the spheroid with $a:=\frac{3}{4}$ starting from $(0, \frac{\pi}{2})$, with the increments $\Delta \varphi_0 = \frac{\pi}{8}$, $t=3$ (on the left) and divided into time segments with $t\in[0, 1] \text{- red}, t\in[1, 2] \text{-blue}, t\in[2, 3] \text{-purple}$.} 
\label{znp_divided}
\end{figure}

\begin{figure}[h]
        \centering
~\includegraphics[width=0.3\textwidth]{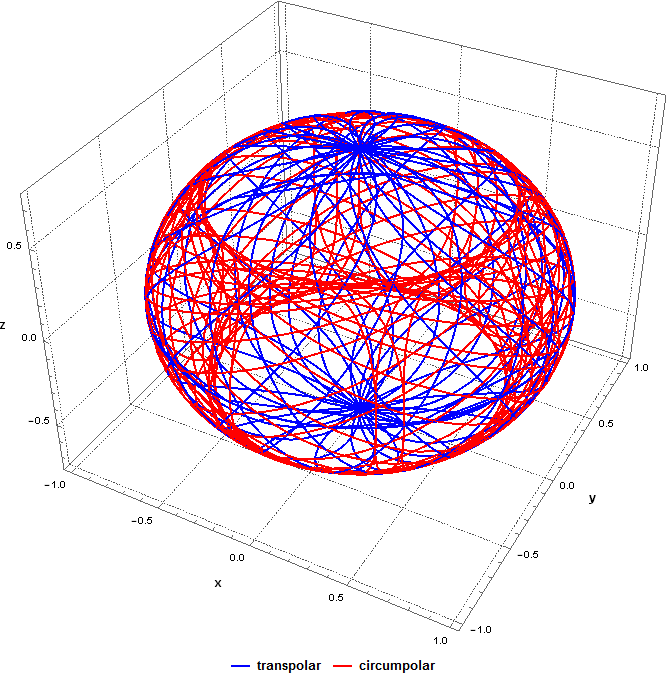}
~\includegraphics[width=0.35\textwidth]{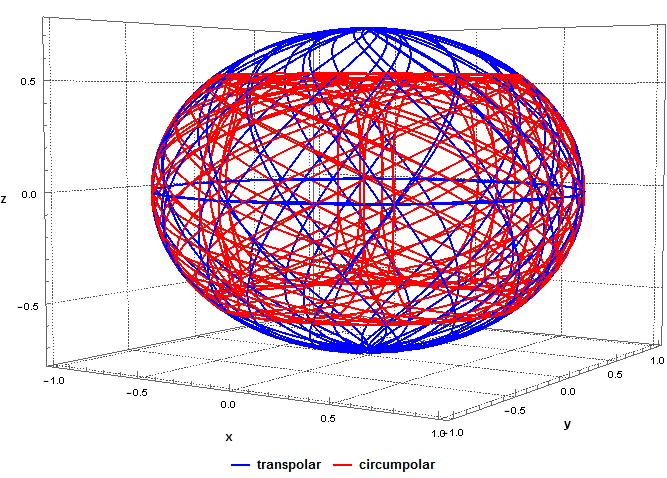}
~\includegraphics[width=0.25\textwidth]{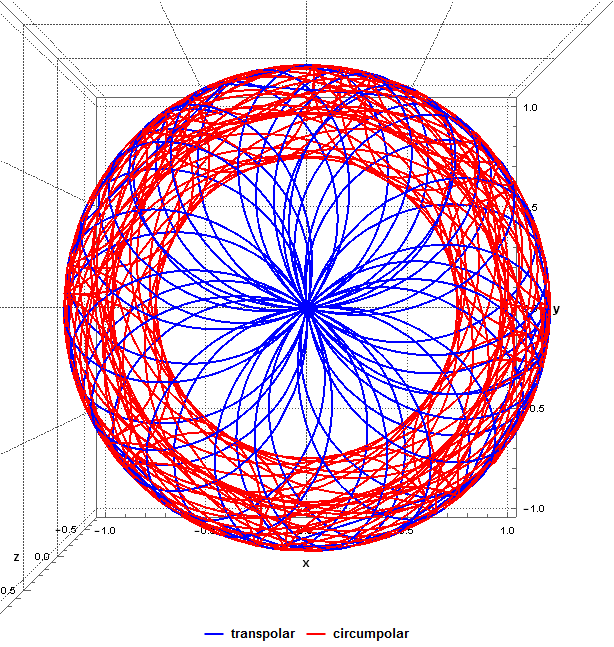}

        \caption{Transpolar (blue) and circumpolar (red) time-optimal paths starting from $(0, \frac{\pi}{2})$ with the increments $\Delta \varphi_0 = \frac{\pi}{4}$, under perturbing infinitesimal rotation with $c:=\frac{5}{7}$, $t=50$.} 
\label{cis_trans}
\end{figure}
\noindent
The form of the initial conditions including the optimal control $\varphi$ under perturbing vector field become $\phi(0)=\phi_0\in[0, 2\pi)$, $\theta(0)=\theta_0\in(0, \pi)$, and for the first derivative 
\begin{equation}
\label{ic_znp}
\dot{\phi}(0)=W^1(\phi_0, \theta_0)+\frac{\cos\varphi_0}{\sin\theta_0}, \quad \dot{\theta}(0)=W^2(\phi_0, \theta_0)-\frac{\sin\varphi_0}{\sqrt{\cos^2\theta_0+a^2\sin^2\theta_0}}. 
\end{equation}
\noindent
where $\varphi=\varphi(t)$ is the angle measured counterclockwise which the vector of the relative velocity forms with a parallel defined by a colatitude $\theta$. The last relations can be derived by direct consideration of the angular equations of motion including the angular representation of the components of ship's own velocity and the background Riemannian metric. 

When the families of the time-optimal paths coming from the same fixed point on the spheroid are considered, $\varphi_0$ plays the role of the parameter which rotates the unit tangent vector of unperturbed Riemannian geodesic. To visualize this we set the increments, for instance $\vartriangle\varphi_0=\frac{\pi}{8}$, in included figures. The time-efficient paths on the spheroid with $a:=\frac{3}{4}$ starting from $(0, \frac{\pi}{2})$, with the increments $\Delta \varphi_0 = \frac{\pi}{8}$, $t=3$ and divided into time segments with $t\in[0, 1] \text{- red}, t\in[1, 2] \text{-blue}, t\in[2, 3] \text{-purple}$ are presented in Figure \ref{znp_divided}. Increasing time we can also observe that the time-optimal paths create transpolar and circumpolar flows on the spheroid. In Figure \ref{cis_trans} we present the transpolar (blue) and circumpolar (red) time-optimal paths starting from $(0, \frac{\pi}{2})$ with the increments $\Delta \varphi_0 = \frac{\pi}{4}$, under perturbing infinitesimal rotation with $c:=\frac{5}{7}$ and $t=50$.

%%%%%%%%%%%%%%%%%%%%%%%%%%%%%%%%%%%%%%%%%%%
\begin{figure}[h]
        \centering
~\includegraphics[width=0.3\textwidth]{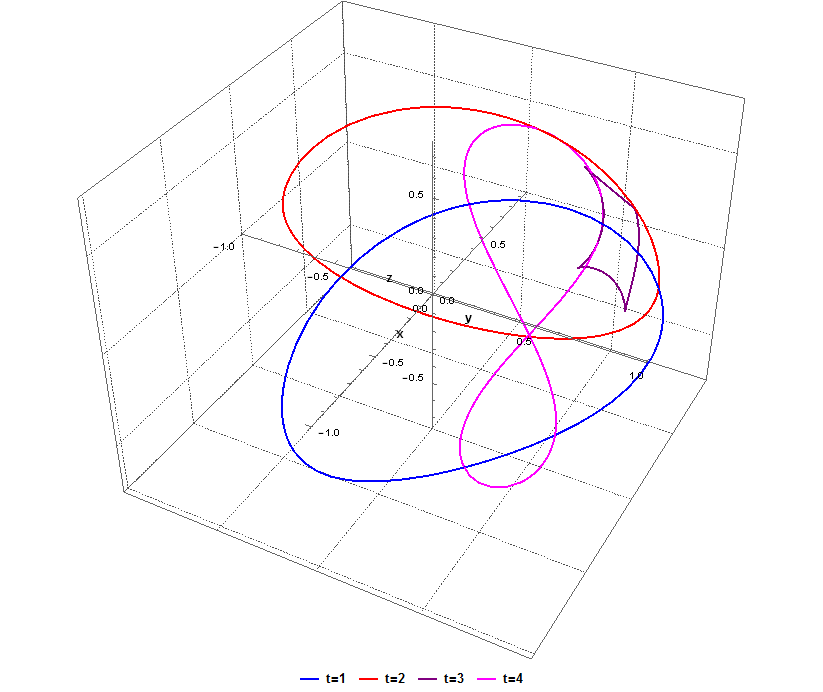}
~\includegraphics[width=0.3\textwidth]{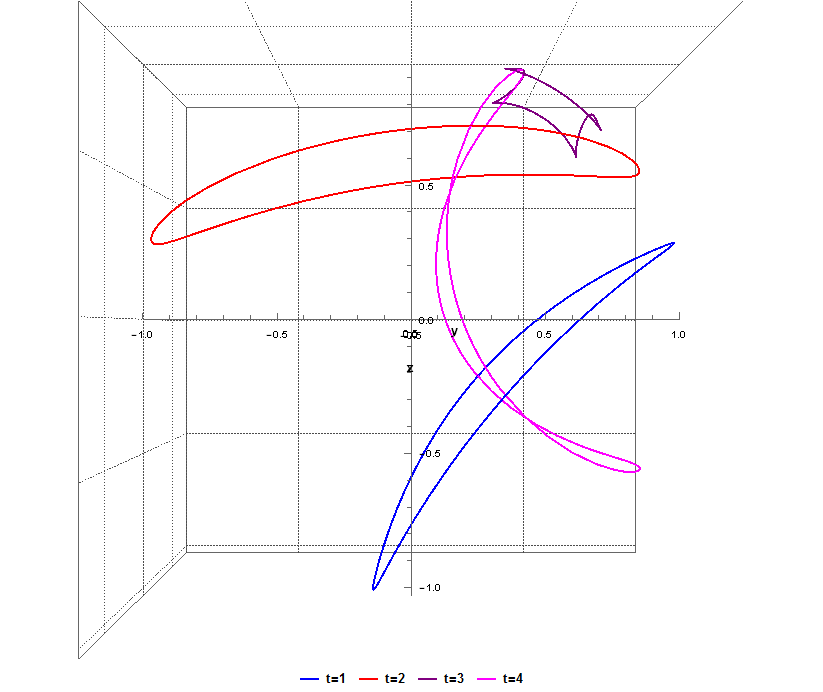}
~\includegraphics[width=0.3\textwidth]{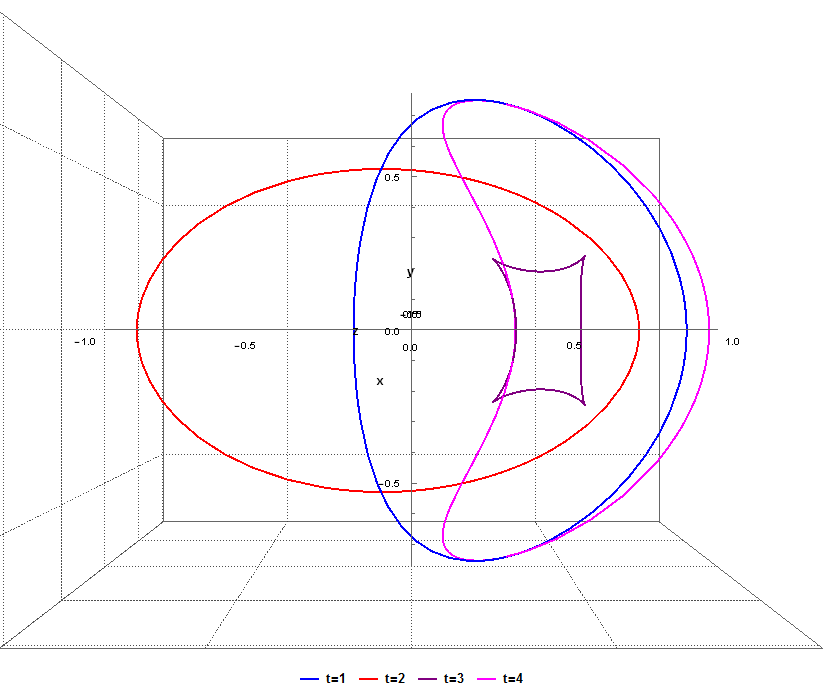}
%~\includegraphics[width=0.35\textwidth]{img/sf_indica_22}

        \caption{The indicatrices under perturbing infinitesimal rotation with $c:=\frac{5}{7}$ for $t=1$ (blue), $t=2$ (red), $t=3$ (purple), $t=4$ (magenta) and the starting point $(0, \frac{\pi}{2})$.}  
\label{indi_1}
\end{figure}
\subsection{Indicatrices and comparison of $h$-, $\alpha$-, $F$- geodesics}

\begin{figure}[h]
        \centering
%~\includegraphics[width=0.3\textwidth]{img/sf_indica_2}
~\includegraphics[width=0.33\textwidth]{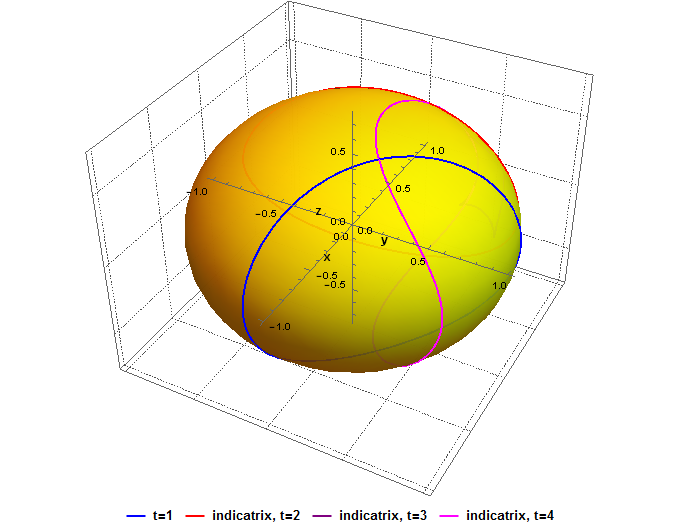}
~\includegraphics[width=0.3\textwidth]{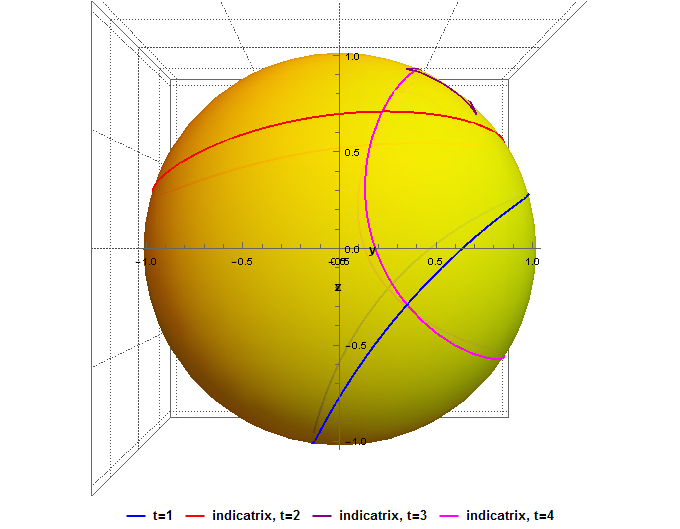}
%~\includegraphics[width=0.45\textwidth]{img/sf_indica_2c}
~\includegraphics[width=0.32\textwidth]{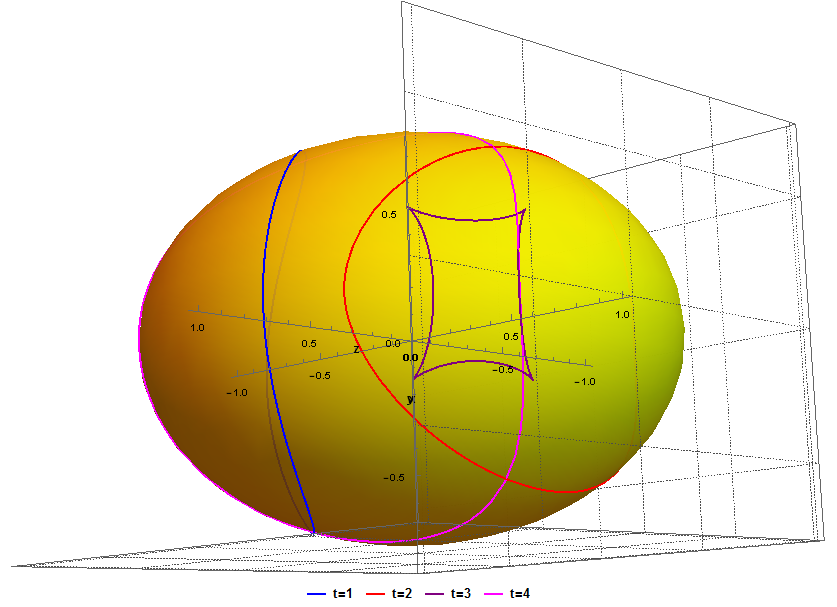}

        \caption{The indicatrices on the spheroid under perturbing infinitesimal rotation with $c:=\frac{5}{7}$ for $t=\{1 \text{(blue)}, 2 \text{(red)}, 3 \text{(purple)}, 4 \text{(magenta)}\}$ and the starting point $(0, \frac{\pi}{2})$.}  
\label{indi_2}
\end{figure}

The unite circle of $h$ in each tangent plane represents the destinations which are reached in one unit of time in the absence of a background wind. This circle is translated rigidly due to the action of the perturbation. Thus, the resulting indicatrices of $F$ are off-centered in comparison to the initial indicatrix and represent the loci of unit time destinations in windy conditions. The indicatrices created by the vector field \eqref{ir} with $c:=\frac{5}{7}$ on the ellipsoid for $t\in\{1, 2, 3, 4\}$ are shown in Figure \ref{indi_1} and Figure \ref{indi_2}. 
%\begin{figure}
   %     \centering
%~\includegraphics[width=0.28\textwidth]{img/sf_indica_2aa}
%~\includegraphics[width=0.32\textwidth]{img/sf_indica_2bb}
%~\includegraphics[width=0.35\textwidth]{img/sf_indica_2cc}
%~\includegraphics[width=0.35\textwidth]{img/sf_indica_2ccc}
   %     \caption{Indicatices for $t=1, 2, 3, 4$} 
%\end{figure}
%\subsection{Flows of the time-optimal paths}
In the presence of a background wind the Riemannian metric $h$ no longer gives the travel time along vectors, but a new metric $F$ on $T\Sigma^2$. Also, in Riemannian geometry two geodesics which pass through a common point in opposite directions necessarily trace the same curve. All reversible Finsler metrics have this property. However, the phenomenon does not extend to nonreversible settings \cite{CR}. 

\begin{figure}[h]
        \centering
~\includegraphics[width=0.45\textwidth]{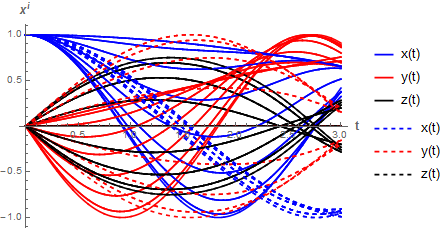}

        \caption{Comparing the solutions $x(t)$-blue, $y(t)$-red, $z(t)$-black before (dashed) and after (solid) perturbation with $\Delta \varphi_0 = \frac{\pi}{8}$, $t=3$ and the starting point $(1, 0, 0)\in\mathbb{R}^3$.} 
\label{geo_znp_xyz}
\end{figure}

\begin{figure}[h]
        \centering
~\includegraphics[width=0.33\textwidth]{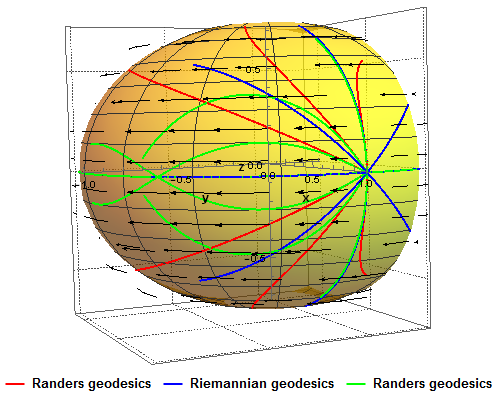}
%~\includegraphics[width=0.35\textwidth]{img/alpha_geo_znp_t=1b}
~\includegraphics[width=0.32\textwidth]{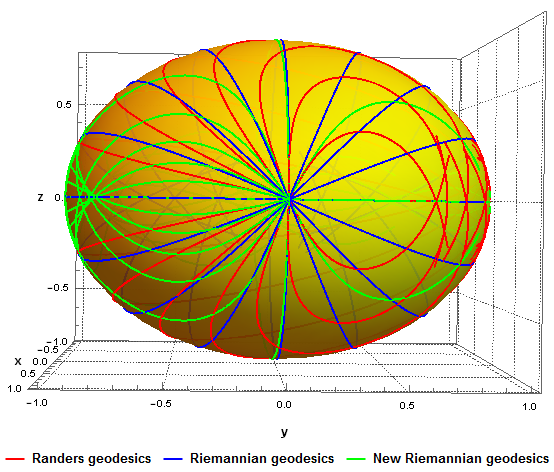}
~\includegraphics[width=0.31\textwidth]{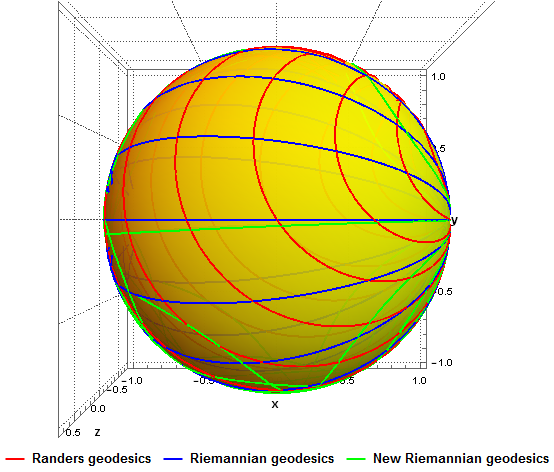}

        \caption{The corresponding background Riemannian (blue), new Riemannian (green) and Randers (red) geodesics starting from $(0, \frac{\pi}{2})$ with $\Delta \varphi_0 = \frac{\pi}{4}$, $t=1$ (on the left ) and with $\Delta \varphi_0 = \frac{\pi}{8}$,  $t=3$.} 
\label{compar_1}
\end{figure}

Much relevant information can be obtained directly from the comparisons of the corresponding geodesics' flows. In Figure \ref{geo_znp_xyz} we can observe the solutions $x(t)$-blue, $y(t)$-red, $z(t)$-black before (dashed) and after (solid) perturbation with $\Delta \varphi_0 = \frac{\pi}{8}$, $t=3$ and the starting point $(1, 0, 0)\in\mathbb{R}^3$. The comparisons of the background Riemannian ($h$-, blue) to the corresponding perturbed new Riemannian ($\alpha$-, green) and the resulting Randers ($F$-, red) geodesics are shown in Figure \ref{compar_1}. One may also check if the background Riemannian geodesic passes or omits the fixed points of the flow of the time-optimal paths (Randers geodesics), for example the ellipsoid's poles.  For a Randers metric $F$ expressed in terms of a Riemannian metric $h$ and a vector field $W$ the relationship between the spray coefficients of $F$, $\alpha$ and $h$ can be found in \cite{colleen_shen, chern_shen}. 
\begin{figure}[h]
        \centering
~\includegraphics[width=0.25\textwidth]{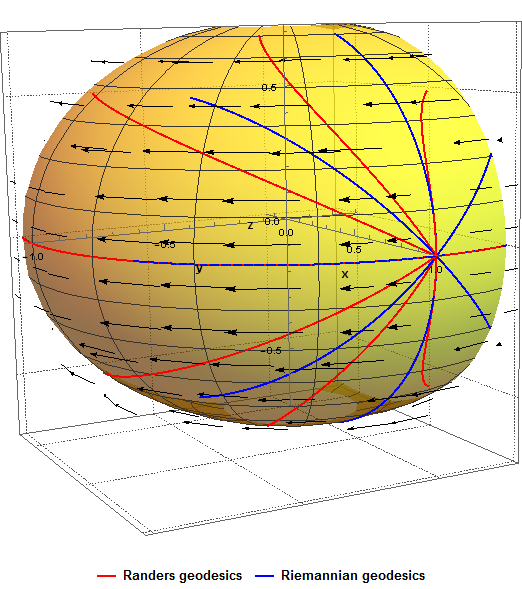}
~\includegraphics[width=0.35\textwidth]{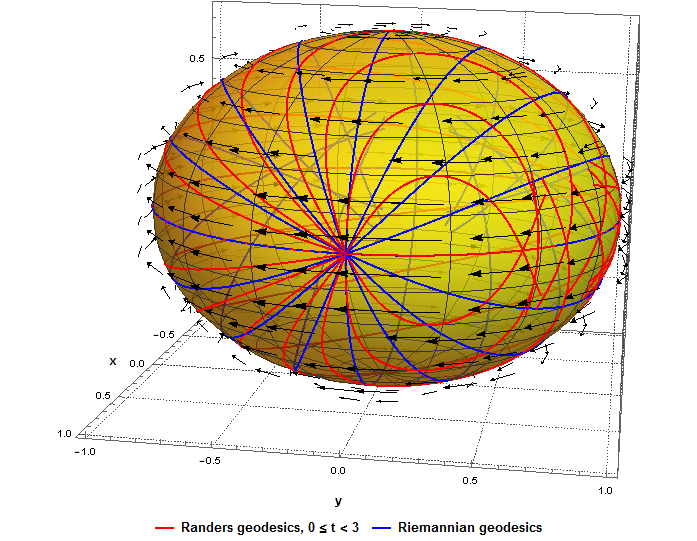}
~\includegraphics[width=0.35\textwidth]{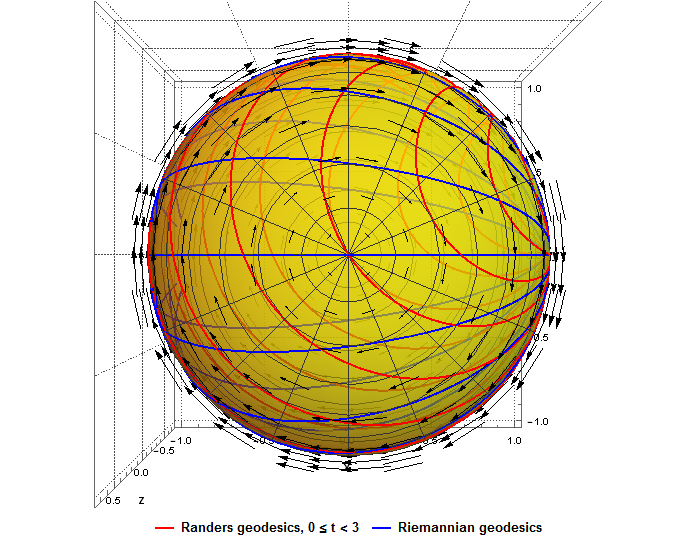}

        \caption{Comparing perturbed (red) and unperturbed (blue) time-optimal paths starting from $(0, \frac{\pi}{2})$, with $\Delta \varphi_0 = \frac{\pi}{4}$, $t=1$ (on the left) and $\Delta \varphi_0 = \frac{\pi}{8}$, $t=3$.} 
\end{figure}
%%%%%%%%%%%%%%%%%%%%%%%%%%%%%%%%%%%%%%%%%%%%%

%\section{Comparing the flows of geodesics}

%%%%%%%%%%%%%%%%%%%%%%%%%%%%%%%%%%%%%%%%%%%

\subsection{Optimal control and the velocities}

\noindent
First, we consider a drift angle $\Psi$  to show a relation between the time-optimal trajectory represented by the resulting course over ground $\Phi$ and its corresponding heading (optimal control) $\varphi$. $\Psi$ determines the angular difference between $\Phi$ and $\varphi$ of the optimal paths. Thus, $\Psi$ shows the effect of acting perturbation on the background Riemannian geodesic. Observe that $|\Psi|<\frac{\pi}{2}$ if $|W|<1$. We apply the convention such that $\Psi$ is positive if the perturbation pushes a navigating ship anticlockwise and negative if it is perturbed clockwise, i.e. $\Psi=\varphi-\Phi$. This corresponds to the real perturbation taken into cosideration in marine navigation like wind or current (stream). Then $\Psi$ is positive if $W$ perturbes a ship to her starboard side and negative if pushes to her port side (the nautical terms for right and left, respectively). However, the angles are taken clockwise from north (a meridian). 

Making use of a dot product with reference to the velocities yields $\cos|\Psi|=|\mathbf{u}|^{-1}|\mathbf{v}|^{-1}(h_{11}u_\phi v_\phi+h_{22}u_\theta v_\theta)$. Recalling that $(\tilde{\phi'}\sin\tilde{\theta})^2+(\tilde{\theta'})^2(\cos^2\tilde{\theta}+a^2\sin^2\tilde{\theta})=1$ and observing that $\tilde{\theta'}=\theta'-W^2$, $\tilde{\phi'}=\phi'-W^1$, where $(\tilde{\phi}, \tilde{\theta})$ state for the solutions to the system of the initial Riemannian geodesics \eqref{sph_geo1} and \eqref{sph_geo2}. Hence, a drift $\Psi \in(-\frac{\pi}{2}, \frac{\pi}{2})$ referring to a time-optimal path on a spheroid $\Sigma^2$ of the semiaxes $(1, 1, a)$, $a>0$, under the influence of a mild perturbation $W(\phi, \theta)$ with $|W(\phi, \theta)|<1$, is given by 
\begin{equation}
\label{eq_drift}
\cos\Psi=\frac{\dot{\phi}(\dot{\phi}-W^1(\phi, \theta))\sin^2\theta+\dot{\theta}(\dot{\theta}-W^2(\phi, \theta))(\cos^2\theta+a^2\sin^2\theta)}{\sqrt{(\dot{\phi}\sin\theta)^2+(\dot{\theta})^2(\cos^2\theta+a^2\sin^2\theta)}},
\end{equation}
\noindent
where $\phi(t), \theta(t)$ determine the geodesics of $F$ defined by \eqref{sph1} and \eqref{sph2}.
In fact the formula \eqref{eq_drift} can be used for any trajectory being a resulting path coming from  any unperturbed path and acting vector field on an arbitrary spheroid, in particular a time-optimal one where the optimality is referred to applied Randers metric $F$. Equivalently, applying the law of cosines in the tangent plane $T\Sigma^2$ yields $|W|^2=|\mathbf{u}|^2+|\mathbf{v}|^2-2|\mathbf{u}||\mathbf{v}|\cos|\Psi|$. One sees immediately that $\cos\Psi=\frac{\lambda+|\mathbf{v}|^2_h}{2|\mathbf{v}|_h}$. Hence, this allows us to prove that   
\begin{equation}
\label{eq_drift2}
\cos\Psi=\frac{1+\sin^2\theta\left[\dot{\phi}^2-(W^1(\phi, \theta))^2\right]+(\cos^2\theta+a^2\sin^2\theta)\left[\dot{\theta}^2-(W^2(\phi, \theta))^2\right]}{2\sqrt{(\dot{\phi}\sin\theta)^2+\dot{\theta}^2(\cos^2\theta+a^2\sin^2\theta)}}.
\end{equation}
%Comparing the right-hand sides of \eqref{eq_drift} and \eqref{eq_drift2} on can show the equivalence 
\noindent
Let us observe that $\Psi$ is positive if $c\dot{\theta}>0$ and negative if $c\dot{\theta}<0$ in the case of perturbing infinitesimal rotation $W(\phi, \theta)=-c\frac{\partial}{\partial\phi}$, $|c|<1$ applied in the example. The result is
\begin{equation}
\Psi=\text{sgn}(c\dot{\theta})\arccos\left (\frac{\dot{\phi}(\dot{\phi}+c)\sin^2\theta+\dot{\theta}^2(\cos^2\theta+a^2\sin^2\theta)}{\sqrt{(\dot{\phi}\sin\theta)^2+\dot{\theta}^2(\cos^2\theta+a^2\sin^2\theta)}}\right).
\label{sgn}
\end{equation}
\noindent
If we allowed stronger windy conditions, namely $|W|\geq1$, then $\Psi\in[-\pi, \pi]$. 
It is therefore of interest to look at $\varphi$. By direct consideration of the linear velocities' triangle in $T\Sigma^2$ and basing on the Randers geodesics of $F$ given by \eqref{sph1} and \eqref{sph2} we find the optimal control (heading) $\varphi$ on perturbed rotational ellipsoid. Optimal control (heading) $\varphi$ on a spheroid $\Sigma^2$ of the semiaxes $(1, 1, a)$, $a>0$ under the influence of an arbitrary mild perturbation $W(\phi, \theta)$ with $|W(\phi, \theta)|<1$, can be given by 
\begin{equation}
\label{opti}
\tan\varphi=\frac{\left(W^2(\phi, \theta)-\dot{\theta}\right)\sqrt{\cos^2\theta+a^2\sin^2\theta}}{\left(\dot{\phi}-W^1(\phi, \theta)\right) \sin\theta}
\end{equation}
Formula \eqref{opti} includes both wind components and excludes the solutions where $|\varphi(t)|=\frac{\pi}{2}$. Actually, having the optimal path in hand which is determined by the Randers geodesics $(\phi, \theta)$, it is sufficient to apply shorter relations which bring the solutions for the entire range $\varphi\in(0, 2\pi)$, namely
\begin{equation}
\label{opti2}
\cos\varphi=\left(\dot{\phi}-W^1(\phi, \theta)\right)\sin\theta, \qquad \sin\varphi=\left(W^2(\phi, \theta)-\dot{\theta}\right)\sqrt{\cos^2\theta+a^2\sin^2\theta},
\end{equation}

\noindent
where $(\phi(t), \theta(t))$ determine the geodesics of $F$ defined by \eqref{sph1} and \eqref{sph2}.  
Computing the control by the inverse function we take into account the quadrant in which the argument lies. It is sufficient to apply \eqref{opti2} with the wind components zeroed to obtain the optimal resulting course over ground $\Phi$. So we thus use the components of the tangent vector to the Randers geodesic which state for the linear components of the resulting velocity. This makes it clear that if $W=0$ then $\Phi=\varphi$ and $\Psi=0$. Both angles are measured counterclockwise from a parallel of latitude determined by the corresponding colatitude $\theta$. The  simple conversion can give the corresponding navigational courses if the angles are taken clockwise from a meridian like an azimuth. 

%COG - The resulting course over ground $\Phi$ represented by the tangent vector to obtained optimal path as follows
%\begin{equation}
%\tan\Phi=\frac{-\dot{\theta}\sqrt{\cos^2\theta+a^2\sin^2\theta}}{\dot{\phi}\sin\theta}
%\end{equation}

%\begin{equation}
%\cos\Phi=\frac{\dot{\phi}\sin\theta}{\sqrt{(\dot{\phi}\sin\theta)^2+(\dot{\theta})^2(\cos^2\theta+a^2\sin^2\theta)}}
%\end{equation}

%\begin{equation}
%\sin\Phi=\frac{-\dot{\theta}\sqrt{\cos^2\theta+a^2\sin^2\theta}}{\sqrt{(\dot{\phi}\sin\theta)^2+(\dot{\theta})^2(\cos^2\theta+a^2\sin^2\theta)}}
%\end{equation}

%%%%%%%%%%%%%%%%%%%%%%%%%%%%%%%%%%%%%%%%

%\subsection{Velocities}

The angular speeds referring to the time-efficient paths are represented by $\omega_{\phi}=\dot{\phi}, \omega_{\theta}=\dot{\theta}$ and the corresponding linear absolute speed's components by $v_{\phi}=\dot{\phi}\sin\theta$, $v_{\theta}=\dot{\theta}\sqrt{\cos^2\theta+a^2\sin^2\theta}$ in the case of unperturbed and perturbed scenario, where the pairs $(\phi, \theta)$ state for the solutions of corresponding Riemannian and Randers geodesic equations, respectively. Comparing the resulting speed 
\begin{equation}
|\mathbf{v}|=\sqrt{\dot{\phi^2}\sin^2\theta+\dot{\theta^2}(\cos^2\theta+a^2\sin^2\theta)}
\label{speed}
\end{equation}
given as the function of time to unit ship's own speed leads to relevant information on when the perturbation acts against or with the navigating ship, increasing or decreasing her resulting speed. Consequently, this influences the total travel time. Substituting \eqref{ic_znp} in \eqref{speed} and rearranging terms, we are thus led to the relation $|\mathbf{v}(\varphi_0)|$. Square of the resulting speed given as the function of the initial control $\varphi_0$, with departing point $(\phi_0, \theta_0)$ in the presence of a perturbation $W$ on the spheroid $\Sigma^2$ reads  
\begin{equation}
%\footnotesize{
\begin{split}
|\mathbf{v}(\varphi_0)|^2=1+(W^1(\phi_0,\theta_0)\sin\theta_0)^2+2W^1(\phi_0,\theta_0)\sin\theta_0\cos\varphi_0 \\ +(W^2(\phi_0,\theta_0))^2(\cos^2\theta_0+(a\sin\theta_0)^2)  -2W^2(\phi_0,\theta_0)\sin\varphi_0\sqrt{\cos^2\theta_0+(a\sin\theta_0)^2}.
\end{split}
%}	
\label{v(fi_0)}
\end{equation}

%\begin{equation}
%\scriptsize{
%v(\varphi_0)=\sqrt{1+(W^1)^2\sin^2\theta_0+2W^1\sin\theta_0\cos\varphi_0+(W^2)^2(\cos^2\theta_0+a^2\sin^2\theta_0)-2W^2\sin\varphi_0\sqrt{\cos^2\theta_0+a^2\sin^2\theta_0}}
%}
%\end{equation}
\noindent
With fixed initial steering one may find the spheroid's parallel $\theta$ for which the speed is extreme or alternatively search for the corresponding heading when commencing the passage from a fixed point on $\Sigma^2$. One sees immediately that for $\tilde{W}(x, y, z)=c(-y, x, z) \rightsquigarrow W(\phi, \theta)=(-c, 0)$ the resulting speed equals
\begin{equation}
|\mathbf{v}(\varphi_0)|=\sqrt{(c\sin\theta_0)^2-2c\sin\theta_0\cos\varphi_0+1}\xrightarrow{(\theta_0=\frac{\pi}{2})}\left\{ \begin{array}{rcl}
1-|c| & \text{(min)} \\ 1+|c| &  \text{(max)}
\end{array}.\right.
\end{equation}

\noindent
The upper and lower limits of $|\mathbf{v}_0|$ differ by at most $2|c|$. In the presented example the initial resulting speed (black) as the function of the initial control angle $\varphi_0\in[0, 2\pi )$ is presented in Figure \ref{v=f(fi_0)} as well as the angular and linear speeds' changes during the ship's passage in Figure \ref{compar_v_liniowe}.

%For the computational convenience, for instance, the two-argument form $\verb"ArcTan[x,y]"$ in the Wolfram programme represents the arc tangent of $\frac{y}{x}$, taking into account the quadrant in which the point $(x, y)$ lies. It therefore gives the angular position expressed in radians of the point measured from the positive $x$ axis.

%For a real number $x$, $arctan(x)$ represents the radian angle measure -[Pi]/2<[Theta]<[Pi]/2 such that tan([Theta])==x.

%%%%%%%%%%%%%%%%%%%%%%%%%%%%%%%%%%%%%%%%

\section{Example with the corresponding individuals}

%\subsection{Time-optimal paths}

We proceed assuming as above for the families of the geodesics that the initial point is determined by $\phi_0=0$, $\theta_0=\frac{\pi}{2}$. By \eqref{sph_geo1} and \eqref{sph_geo2} and setting $a:=\frac{3}{4}$ the background Riemannian geodesic equations of $\Sigma^2$ become 
\begin{equation}
\ddot{\phi}+2 \dot{\theta} \dot{\phi} \cot\theta=0,  \qquad \ddot{\theta}-\frac{\left(7 \dot{\theta} ^2+16 \dot{\phi} ^2\right)\sin 2 \theta}{7 \cos 2 \theta +25}=0.
\label{sph_geo1_ex}
\end{equation}
\ 
The initial conditions for the individual are complemented by $\dot{\phi} =\frac{1}{2}, \dot{\theta}(0)=-\frac{2}{\sqrt{3}}$, so $|(\frac{1}{2}, -\frac{2}{\sqrt{3}})|=1$. Thus, a ship begins with the initial heading $\varphi_0=\Phi_0=\frac{\pi}{3}$ with unit speed through the water. 
\begin{figure}[h]
        \centering
~\includegraphics[width=0.32\textwidth]{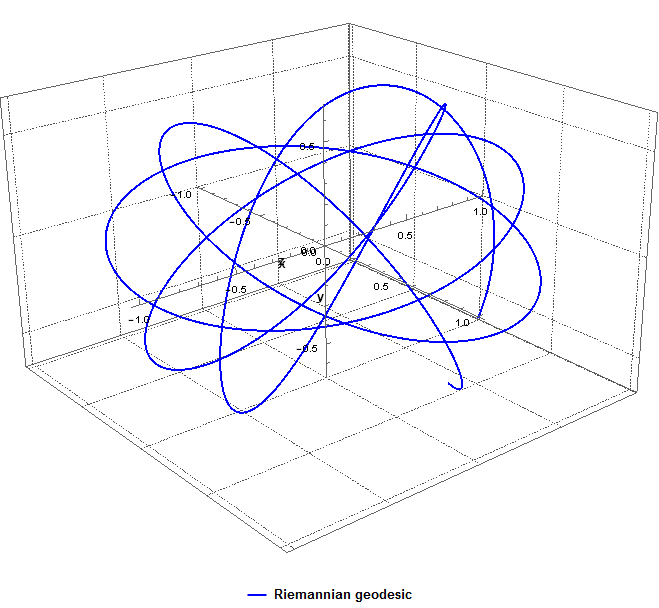}
~\includegraphics[width=0.33\textwidth]{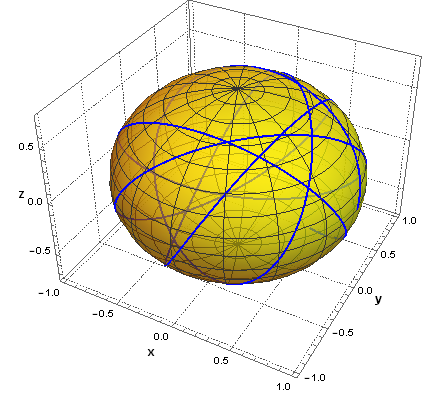}
%~\includegraphics[width=0.3\textwidth]{img/sf_geo3a}
~\includegraphics[width=0.27\textwidth]{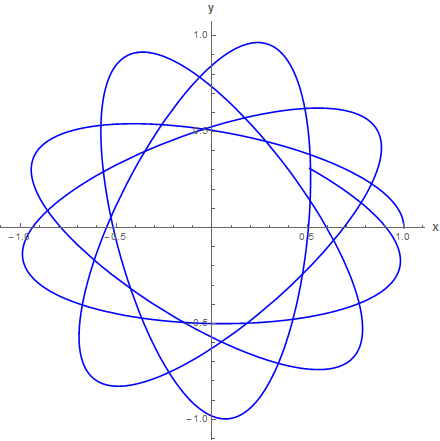}

        \caption{The Riemannian geodesic (the individual) departing from $(0, \frac{\pi}{2})$ and its planar $xy$-projection, with $t=25$.} 
\label{sf_geo}
\end{figure}
This corresponds to the heading (true course) 030$^\circ$ $\approx$ NEbN (northeast by north) in the rhumb system with respect to true north on a compass rose applied to the oblate spheroid. 
\noindent
The Riemannian geodesic departing from $(0, \frac{\pi}{2})$ and its planar $xy$-projection, with $t=25$ are presented in Figure \ref{sf_geo}. 
\noindent
Plugging also the constant $c$  and rearranging terms we are thus led to the relation for $L\text{=}\frac{1}{2}F^2$
%\begin{equation}
%\footnotesize{
%F(\phi ,\theta; u,v)=\frac{7 \left(\sqrt{-64 \cos (2 \theta ) \left(98 u^2-71 v^2\right)+6272 u^2+350 v^2 \cos (4 \theta )+7650 v^2}+80 u \sin ^2(\theta )\right)}{8 (25 \cos (2 \theta )+73)}.
%}
%\end{equation}
\begin{equation}
%\footnotesize{
L(\phi ,\theta ,u,v)=\frac{49 \left(\sqrt{-4 \cos 2 \theta \left(98 u^2-71 v^2\right)+392 u^2+\frac{175}{8} v^2 \cos 4 \theta+\frac{3825 v^2}{8}}+20 u \sin ^2\theta \right)^2}{8 (25 \cos 2 \theta +73)^2}.
%}
\end{equation}
\noindent
Next, by \eqref{alpha_1ex} and \eqref{alpha_2ex} additionally we obtain the geodesic equations for the new Riemannian metric $\alpha$ in the following form
\begin{equation}
\ddot{\phi}+\frac{2 \dot{\phi} \dot{\theta} \cot \theta  \left(49 \csc ^2\theta +25\right)}{49 \csc ^2\theta -25}=0,  
\end{equation}
\begin{equation}
\ddot{\theta}+\frac{2 \sin 2 \theta  \left(25 \cos 2 \theta  \left(784 \dot{\phi}^2+57 \dot{\theta}^2\right)-96432 \dot{\phi}^2+4161 \dot{\theta}^2\right)}{(7 \cos 2 \theta +25) (25 \cos 2 \theta +73)^2}=0.
\end{equation}

\ 

\begin{figure}[h]
        \centering
~\includegraphics[width=0.34\textwidth]{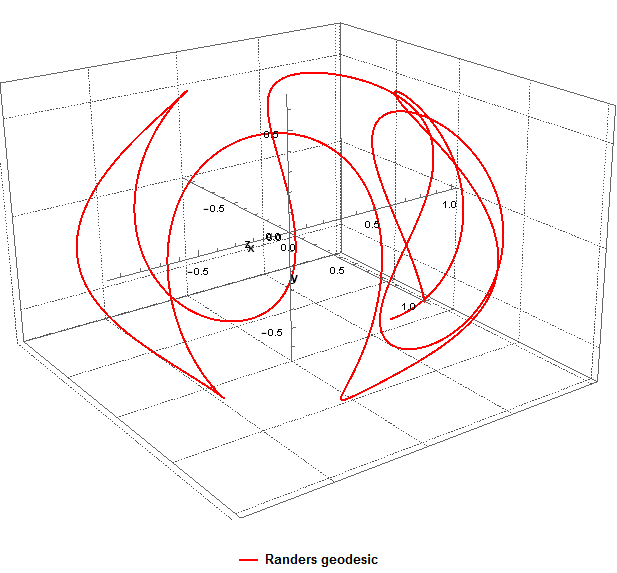}
%~\includegraphics[width=0.3\textwidth]{img/sf_znp_t=25_3}
~\includegraphics[width=0.31\textwidth]{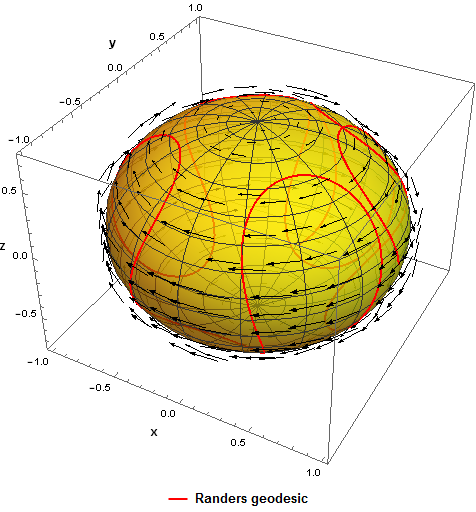}
%~\includegraphics[width=0.31\textwidth]{img/sf_znp_t=25_4}
~\includegraphics[width=0.3\textwidth]{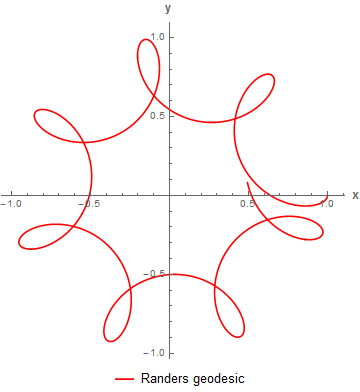}

        \caption{The time-optimal path (the individual) starting from $(0, \frac{\pi}{2})$ under rotational perturbation with $c:=\frac{5}{7}$, $t=25$ and its planar $xy$-projection with $t=20$.} 
\label{sf_znp}
\end{figure}
\noindent
Now, recalling $1-$form $\beta$ by \eqref{sph2} we thus get the final equations of the time-optimal paths. 
%\begin{equation}
%\footnotesize{
%\begin{split}
%G(\phi ,\theta; u,v)=v \cot \theta  \left(\sqrt{-64 \cos 2 \theta  \left(98 u^2-71 v^2\right)+6272 u^2+350 v^2 \cos 4 \theta +7650 v^2}-40 u \cos 2 \theta +40 u\right)^3 \\
%\left(5 \sqrt{-64 \cos 2 \theta  \left(98 u^2-71 v^2\right)+6272 u^2+350 v^2 \cos 4 \theta +7650 v^2}+784 u\right)
%\end{split}
%}
%\end{equation}
%\begin{equation}
%\tiny{
%H(\phi ,\theta; u,v)=16 (25 \cos 2 \theta +73) \left(-8000 u^3 \cos 6 \theta +644480 u^3+6736 u^2 \sqrt{-64 \cos 2 \theta  \left(98 u^2-71 v^2\right)+6272 u^2+350 v^2 \cos 4 \theta +7650 v^2}+3825 v^2 \sqrt{-64 \cos 2 \theta  \left(98 u^2-71 v^2\right)+6272 u^2+350 v^2 \cos 4 \theta +7650 v^2}+5 \cos 4 \theta  \left(47232 u^3+240 u^2 \sqrt{-64 \cos 2 \theta  \left(98 u^2-71 v^2\right)+6272 u^2+350 v^2 \cos 4 \theta +7650 v^2}+35 v^2 \sqrt{-64 \cos (2 \theta ) \left(98 u^2-71 v^2\right)+6272 u^2+350 v^2 \cos 4 \theta +7650 v^2}-23064 u v^2\right)-4 \cos 2 \theta  \left(218160 u^3+1984 u^2 \sqrt{-64 \cos 2 \theta  \left(98 u^2-71 v^2\right)+6272 u^2+350 v^2 \cos 4 \theta +7650 v^2}-568 v^2 \sqrt{-64 \cos 2 \theta  \left(98 u^2-71 v^2\right)+6272 u^2+350 v^2 \cos 4 \theta +7650 v^2}+49215 u v^2\right)-10500 u v^2 \cos 6 \theta +322680 u v^2\right)
%}
%\end{equation}
Let
\begin{equation*}
\tilde{\tau} = \sqrt{\dot{\theta} ^2 (4544 \cos 2 \theta +350 \cos 4 \theta +7650)+12544 \dot{\phi} ^2\sin ^2\theta  }.
\end{equation*}
%%koniec zmiennych
Hence, 
\begin{equation}
\ddot{\phi} +\frac{\left(5 \tilde{\tau} +784 \dot{\phi} \right)\dot{\theta}  \cot \theta  \csc ^2\theta}{8 \left(49 \csc ^2\theta -25\right)}=0,
\end{equation}

\begin{equation}
\ddot{\theta} -
\frac{
\splitfrac{ 
	 \left(\tilde{\tau}+80  \dot{\phi}\sin ^2\theta \right)^3 
	\left(392 \dot{\phi} \left(5 \tilde{\tau}+(492-100 \cos 2 \theta ) \dot{\phi} \right)
	\right.}{\left. + \dot{\theta} ^2 (53950 \cos 2 \theta )+4375 \cos 4 \theta +87303)\right)\sin 2 \theta 
}
}{
\splitfrac{ 
2  (7 \cos 2 \theta +25) (25 \cos 2 \theta +73)^2\left(\dot{\theta} ^2 (2272 \cos 2 \theta +175 \cos 4 \theta +3825)  \right.}{ \left. \left(\tilde{\tau}+240 \dot{\phi} \sin ^2\theta  \right)
+64 \dot{\phi} ^2\sin ^2\theta   \left((173-75 \cos 2 \theta ) \tilde{\tau}-80 \sin ^2\theta  (25 \cos 2 \theta -319) \dot{\phi} \right)\right)
}
}=0.
\end{equation}

\begin{figure}[h]
        \centering
~\includegraphics[width=0.3\textwidth]{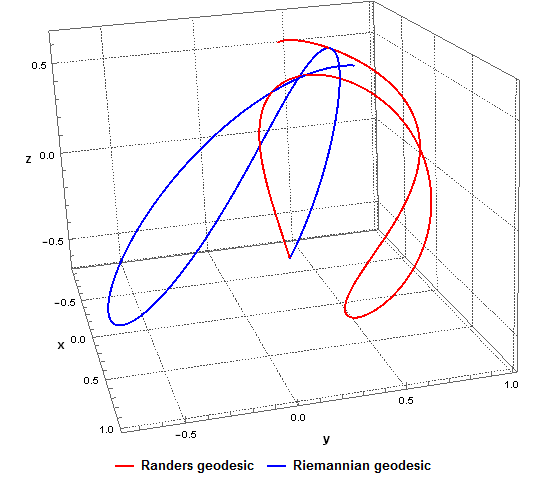}
~\includegraphics[width=0.3\textwidth]{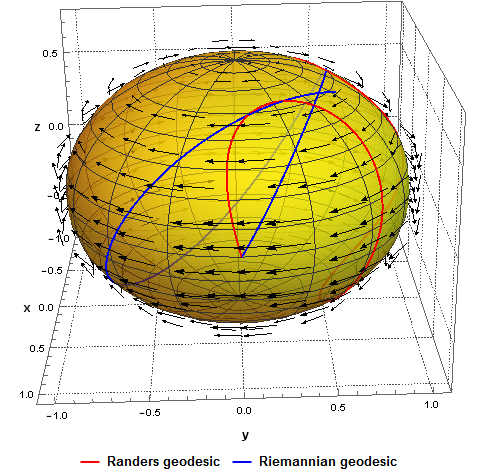}
%~\includegraphics[width=0.3\textwidth]{img/compar_geo_znp_t=7a}
~\includegraphics[width=0.3\textwidth]{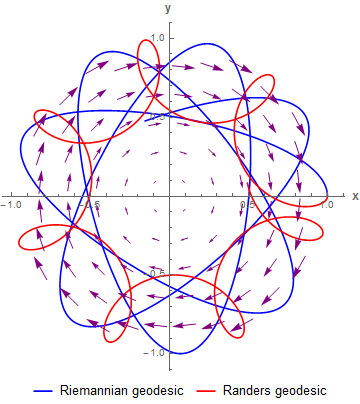}
        \caption{The background Riemannian geodesic vs. the corresponding Randers geodesic with $t=7$ and their $xy$-projections with $t=20$.} 
\label{compar_geo_znp_t=7_2}
\end{figure}

\begin{figure}[h]
        \centering
~\includegraphics[width=0.36\textwidth]{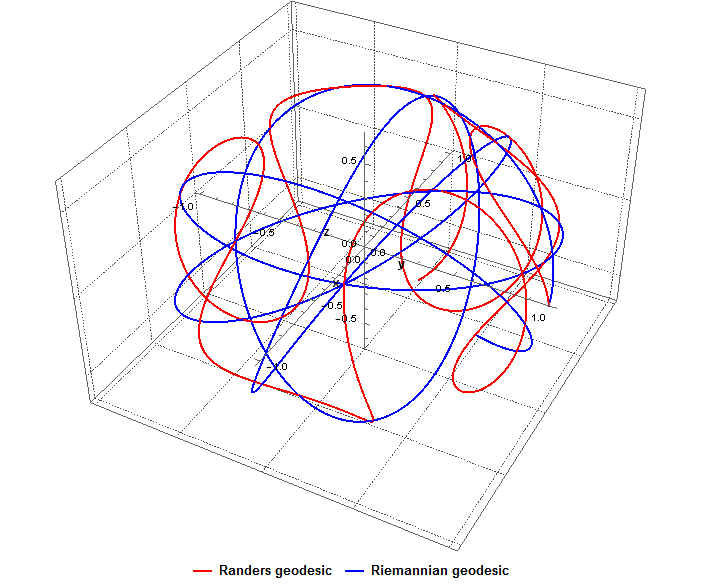}
%~\includegraphics[width=0.33\textwidth]{img/compar_geo_znp_t=25_4}
~\includegraphics[width=0.3\textwidth]{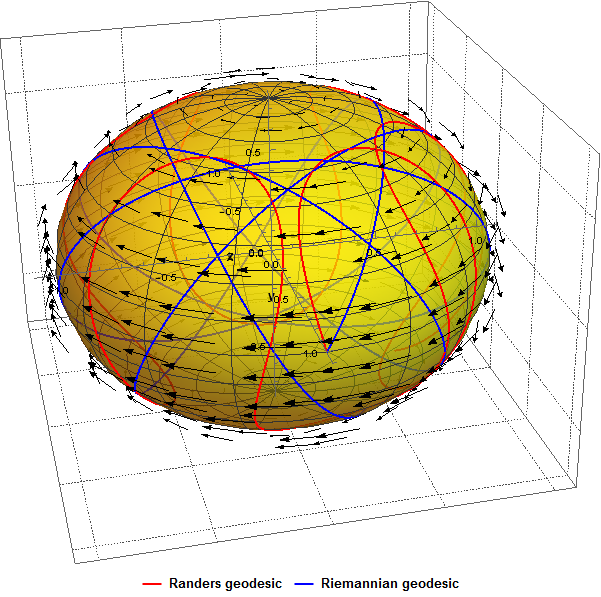}
~\includegraphics[width=0.3\textwidth]{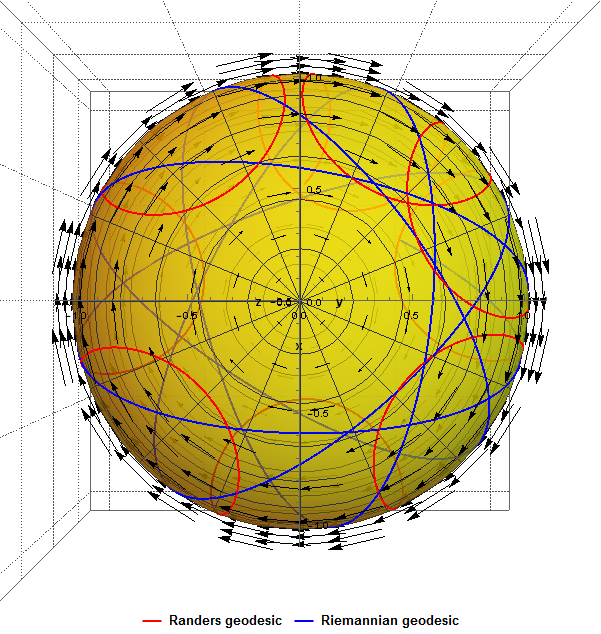}

        \caption{Comparing the individuals before (blue) and after (red) infinitesimal rotation with $c:=\frac{5}{7}$, $t=25$.} 
\label{compar_geo_znp_t=25}
\end{figure}

\begin{figure}[h]
        \centering
~\includegraphics[width=0.47\textwidth]{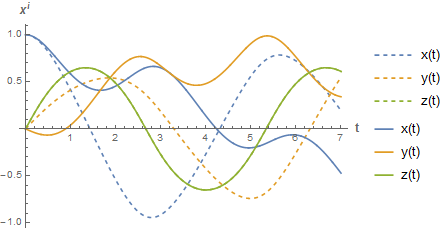}
%~\includegraphics[width=0.40\textwidth]{img/compar_rozw}
~\includegraphics[width=0.32\textwidth]{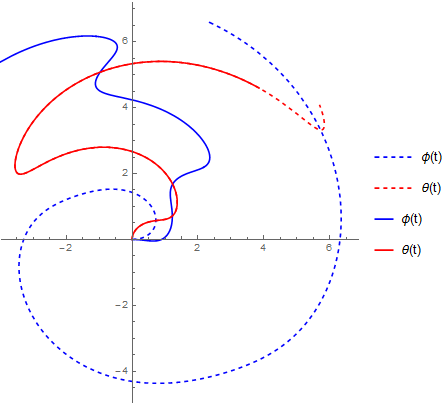}

        \caption{The solutions before (dashed) and after (solid) perturbation: in the base $(x, y, z)$ (on the left) and the polar plot of the solutions in the base $(\phi, \theta)$, $t=7$.} 
\label{compar_rozw_xyz}
\end{figure}

\noindent
To solve the system of the Randers geodesics' equations the initial conditions are set according to \eqref{ic_znp}, i.e. $\phi_0=0$, $\theta_0=\frac{\pi}{2}$, $\dot{\phi}(0)=-\frac{3}{14}, \dot{\theta}(0)=-\frac{2}{\sqrt{3}}$. The resulting geodesic creates the circumpolar time-efficient path on the ellipsoid since $0<\theta(t)<\pi$ for any $t$. We can obtain much relevant information directly from the graphs' comparisons. We thus present the graph of the background Riemannian geodesic and the corresponding Randers geodesic with $t=7$ and their $xy$-projections, with $t=20$, which are compared in Figure \ref{compar_geo_znp_t=7_2} and both paths, with $t=25$, in Figure \ref{compar_geo_znp_t=25}.
The solutions before (dashed) and after (solid) perturbation: in the base $(x, y, z)$ and the polar plot of the solutions in the base $(\phi, \theta)$, $t=7$ are shown in Figure \ref{compar_rozw_xyz}. 

\begin{figure}[h!]
        \centering
~\includegraphics[width=0.5\textwidth]{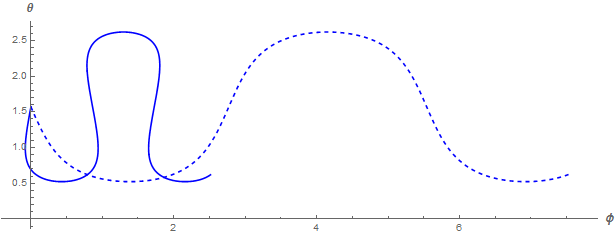}
~\includegraphics[width=0.19\textwidth]{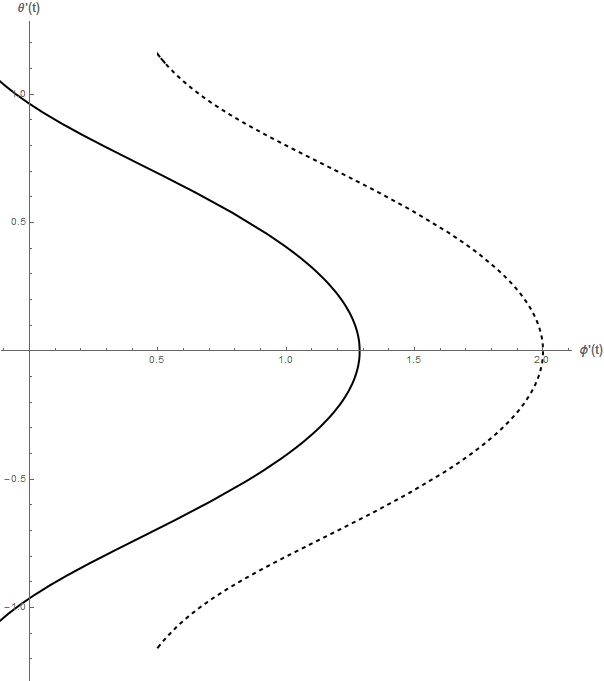}
~\includegraphics[width=0.24\textwidth]{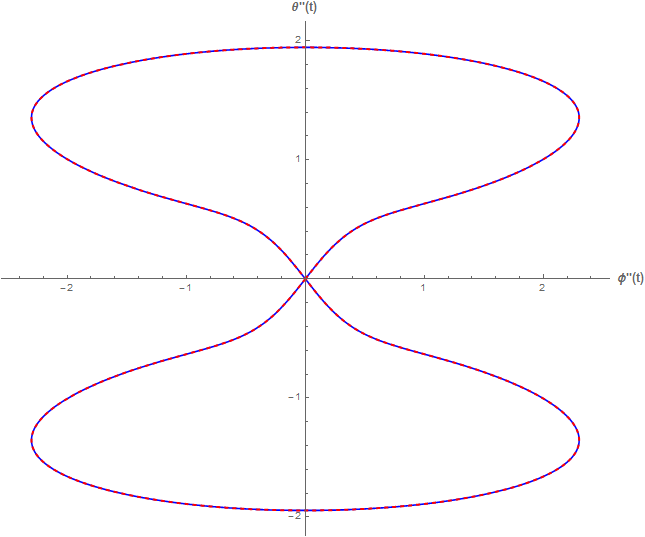}

        \caption{The parametric plots of the solutions ($t=7$,on the left), their first derivatives ($t=3$, in the middle) and second derivatives ($t=5.35$,on the right), before (dashed) and after (solid) perturbation. } 
\label{compar_fi_teta_param}
\end{figure}

We can observe that $z(t)$ and $\theta(t)$ did not change due to acting perturbation. The same holds for the meridian components of the linear velocities $v_\theta$ before and after perturbation. This is the consequence of the fact that $\tilde{W}^3(x, y, z)=0$ what implies $W^2(\phi, \theta)= 0$. The effect of perturbing wind $W(\phi, \theta)=-\frac{5}{7}\frac{\partial}{\partial\phi}$ causes that the time-optimal path in windy conditions is modified anticlockwise at the begining as shown in Figure \ref{sf_znp}. 
Additionally, we also present the parametric plots of the solutions ($t=7$, on the left), their first derivatives ($t=3$, in the middle) and second derivatives ($t=5.35$, on the right), before (dashed) and after (solid) perturbation in Figure \ref{compar_fi_teta_param}. The perturbation causes that the resulting course $\Phi_0\approx 103.9^{\circ}$ at the departure while the optimal control $\varphi_0=\frac{\pi}{3}$. Hence, the initial drift is negative, $\Psi_0\approx-43.9^\circ$.

\begin{figure}[h]
        \centering
~\includegraphics[width=0.58\textwidth]{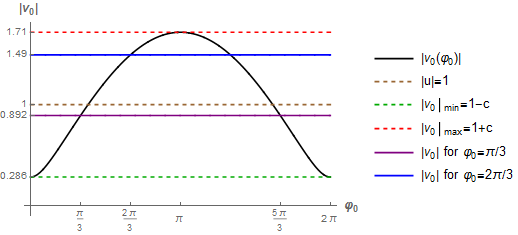}
%UWAGA: TU ZMIENIŁEM NAZWE OBRAZKA GDYZ ZAMIENIAŁ NAWIASY OKRAGLE NA DOLNE SPACJE W EDYTORZE ARXIV
%~\includegraphics[width=0.55\textwidth]{img/v=f(fi_0)_2}

        \caption{The initial resulting speed $|\mathbf{v}_0|$ (black) as the function of the initial control angle $\varphi_0\in[0, 2\pi )$.} 
\label{v=f(fi_0)}
\end{figure}

\begin{figure}[h]
        \centering
~\includegraphics[width=0.49\textwidth]{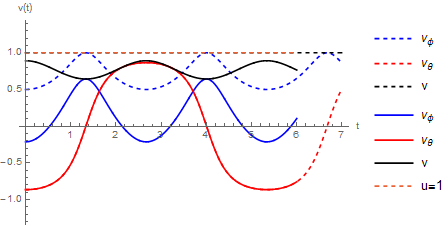}
~\includegraphics[width=0.46\textwidth]{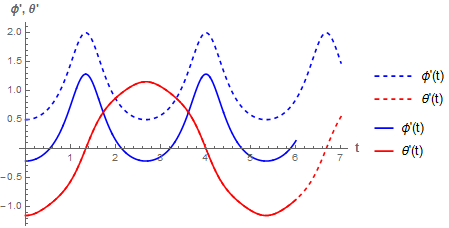}

        \caption{The linear (on the left) and angular speeds as the functions of time, before (dashed) and after (solid) perturbation; the resulting linear speed is shown in black.} 
\label{compar_v_liniowe}
\end{figure}

Having applied \eqref{v(fi_0)} we thus get the graph of the initial resulting speed $|\mathbf{v}_0|$ (black) as the function of the initial control angle $\varphi_0\in[0, 2\pi )$ what is presented in Figure \ref{v=f(fi_0)}. The upper and lower limits of $|\mathbf{v}_0|$ differ by at most $2|c|=\frac{10}{7}$. 
Due to the fact that a ship commences against the perturbation the speed must be decreased. The graph of the linear and angular speeds before (dashed) and after (solid) perturbation are presented in Figure \ref{compar_v_liniowe}. This shows the translation of the graphs of the angular speeds $\dot{\phi}(t)$ by the constant vector $[-\frac{5}{7}, 0]$ as the result of $W^1(\phi, \theta)=-c=const.$ for any $t$. The resulting speed equals $|\mathbf{v}_0|=\frac{\sqrt{39}}{7}\approx 0.892$ at the departure. It is also lower than ship's own speed through the water during the entire voyage, namely $|\mathbf{v}(t)|<1$ for any $t$. For comparison, if we left with the current, changing the sign of $u$, i.e. $\dot{\phi}(0)=-\frac{1}{2}$, then the optimal heading  would be $\varphi_0=\frac{2\pi}{3}$ what corresponds to 330$^\circ$ $\approx$ NWbN (northwest by north) in the rhumb system with respect to true north on the spheroid's compass rose. Then the tangent vector to the time-optimal path $(u, v):=(-\frac{17}{14},-\frac{2}{\sqrt{3}})$ what causes that $\Phi_0\approx 144.5^{\circ}$ over ground and $\Psi_0\approx-24.5^\circ$. Consequently, by \eqref{v(fi_0)} the resulting speed increases, i.e. $|\mathbf{v}_0|=\frac{\sqrt{109}}{7} \approx 1.491$. 
\begin{figure}[h!]
        \centering
~\includegraphics[width=0.53\textwidth]{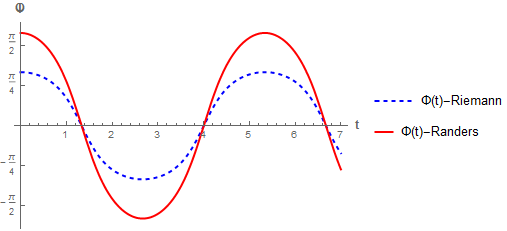}
%~\includegraphics[width=0.2\textwidth]{img/compar_cog2}
~\includegraphics[width=0.33\textwidth]{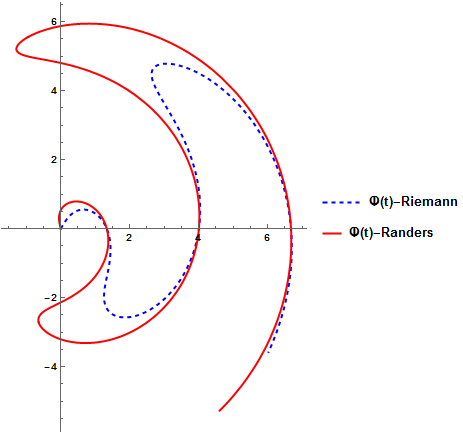}
%~\includegraphics[width=0.35\textwidth]{img/compar_cog1}

        \caption{The time-optimal course over ground $\Phi$ under windy conditions, before (blue dashed) and after (red solid) perturbation in the Cartesian graph and the corresponding polar plot, $t=7$.} 
\label{compar_cog_2}
\end{figure}

\begin{figure}[h!]
        \centering
%~\includegraphics[width=0.45\textwidth]{img/compar_znp_cog_tc2}
%~\includegraphics[width=0.2\textwidth]{img/compar_tc2}
%~\includegraphics[width=0.35\textwidth]{img/compar_tc1}
~\includegraphics[width=0.55\textwidth]{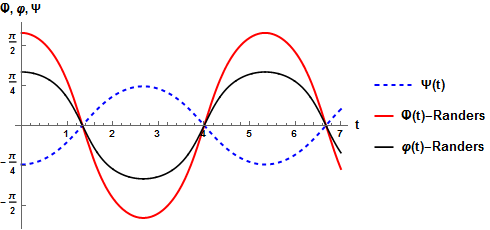}

       % \caption{Optimal "course over ground" $\Phi$ (red) and "heading" (black) in perturbed scenario.} 
        \caption{The drift angle $\Psi(t)$ (blue dashed), the optimal control $\varphi(t)$ (black) and the optimal resulting course over ground $\Phi(t)$ (red) in the perturbed scenario.} 
\label{drejfus}
\end{figure}

In the example the graphs of $\varphi$ as the functions of time are the same for the scenarios before and after perturbation while the investigations of $\Phi(t)$ differ. This means that we observe the continuous nonzero wind effect in almost entire passage (excluding $t$ for which $\varphi(t)=\Phi(t)$) and we let the perturbation to deflect our route (passive approach). This corollary also arises clearly from \eqref{sgn} applied to the example. A different task (active approach) would be if a ship reacts continously to neutralize the effect of the perturbation and adjusts $\varphi(t)$ such that she follows a preset optimal track over ground determined by $\Phi(t)$. The time-optimal course over ground $\Phi$ before and after perturbation in the Cartesian graph and the corresponding polar plot are presented in Figure \ref{compar_cog_2}. These curves arise from \eqref{opti2}. The control (steering course) $\varphi(t)$ according to \eqref{opti2} guarantees that the ship proceeds the time-optimal path on $\Sigma^2$. Finally, by \eqref{sgn} we thus obtain the graph the drift angle $\Psi(t)$ which can be compared to the optimal control $\varphi(t)$ and the optimal resulting course over ground $\Phi(t)$  in the perturbed scenario what is shown in Figure \ref{drejfus}. $|\Psi(t)|\in[0, \approx 43.9^\circ]$ so our ship is drifted to her port side (counterclockwise) and starboard side (clockwise) alternately during her voyage on the perturbed ellipsoidal sea.   

%NA WYK PREDKOSCI KATOWYCH WIDAC STALE PRZESINIECIE BO WIATR KATOW JETS STALY

%\begin{figure}
   %     \centering
%~\includegraphics[width=0.3\textwidth]{img/sf_geo}
%~\includegraphics[width=0.3\textwidth]{img/sf_geo2}
 %       \caption{Riemannian geodesic - individual, xy-projection, t=7} 
%\end{figure}

%%%%%%%%%%%%%%%%%%%%%%%

%\subsection{Optimal control $\varphi$ and the resulting "courses over ground" $\Phi$}

%\begin{figure}
%        \centering
%~\includegraphics[width=0.55\textwidth]{img/drift_2}
%~\includegraphics[width=0.55\textwidth]{img/drift_1}
%~\includegraphics[width=0.55\textwidth]{img/drift_3}
   %     \caption{Drift $\Psi(t)$ (blue dashed), optimal control $\varphi(t)$ (black) and optimal resulting course $\Phi(t)$ in perturbed scenario.} 
%\end{figure}

%%%%%%%%%%%%%%%%%%%%%%%%%%%%%%%%%%%%%%%%%%%

\section{Conclusions}

The analysis including the optimal control and drift angle can be applied to an arbitrary oblate ($0<a<1$) or prolate ($a>1$) ellipsoid in the presence of acting mild vector field $W$, $|W|<1$. The solutions require the restriction of a strong convexity which arises from the assumption of applied Proposition 1.1 of \cite{colleen_shen}. The theorem establishes the direct relation between the Randers geodesics and the time-optimal paths as the solutions to the navigation problem. Solving the system of the Randers geodesics of $F$ given by \eqref{sph1} and \eqref{sph2} and applying the formula for $\varphi$ yield the time-optimal control on perturbed spheroid which is widely applied as the geometric model in real applications mentioned above. In that way it also gives rise to obtain the spheroidal analogue of the classical formula of Zermelo for the optimal heading which implicitly solved the problem of finding the shortest time paths with the Euclidean background, namely in $\mathbb{R}^2$ and $\mathbb{R}^3$  (cf. \cite{zermelo2, zermelo, caratheodory}). 
In order to complement the solution to the problem it will be reasonable to consider also the perturbation which satisfies $|W|_h\geq1$. The obstacle which appears is to obtain the explicit solutions since the computational analyses involving Randers spaces are generally difficult. However the numerical ones followed by the relevant schemes give useful information on the geometric properties of the time-efficients paths including the optimal control what is of our research interest. Regarding navigational applications it would be more convenient to adopt the orientation of $\varphi$ such that it is represented by the angle with respect to north determined by the ellipsoid's meridians or to convert the obtained $\varphi$ to the corresponding azimuth. By applying Clairaut's theorem to the spheroid  $\Sigma^2$ one can also deduce the azimuths of the background unperturbed geodesics. We conclude by noting that a number of non-Riemanniam Randers metrics that are either Einstein or Ricci-constant which include other  surfaces of revolution as the examples with Riemannian-Einstein navigation data $(h, W)$ can be found in Section 4 of \cite{bao_cr}. 

%\bigskip

%\ 

%\noindent
%\textbf{Acknowledgement} \quad The research was supported by a grant from the Polish National Science Center under research project number 2013/09/N/ST10/02537.  

\bigskip

%\noindent
%\textbf{Conflict of Interest:} The author declares that there is no conflict of interest.

\bigskip

%\noindent
%Piotr Kopacz\\
%Faculty of Mathematics and Computer Science, Jagiellonian University,\\
%6, Prof. S. Łojasiewicza, 30 - 348, Kraków, Poland\\
%and \\
%Faculty of Navigation, Gdynia Maritime University,\\
%3,  Al. Jana Pawla II, 81-345, Gdynia, Poland.\\
%E-mail: \verb"piotr.kopacz@im.uj.edu.pl"\\

%\noindent
%\textsc{Faculty of Mathematics and Computer Science, Jagiellonian University,\\
%6, Prof. S. Łojasiewicza, 30 - 348, Kraków, Poland}\\
%and \\
%\textsc{Faculty of Navigation, Gdynia Maritime University,\\
%3,  Al. Jana Pawla II, 81-345, Gdynia, Poland}\\

\noindent
%\textit{E-mail address:} \texttt{piotr.kopacz@im.uj.edu.pl}

%\bibliographystyle{apalike}
%\bibliographystyle{model1b-num-names}
%\bibliographystyle{elsarticle-num}
\bibliographystyle{siam}
\bibliography{pk_bjga}

\end{document}